\documentclass[12pt]{amsart}
\usepackage{latexsym, amsmath, amscd, amssymb, amsthm, bm}
\usepackage{mathrsfs, euscript}
\usepackage{epsfig}
\usepackage{psfrag}
\usepackage[all]{xypic}
\newcommand{\scr}[1]{\bm{\EuScript{#1}}}
\DeclareMathOperator{\supp}{supp}
\DeclareMathOperator{\reg}{reg}
\DeclareMathOperator{\pos}{pos}
\DeclareMathOperator{\intt}{int}
\DeclareMathOperator{\Hilb}{Hilb}
\DeclareMathOperator{\Spec}{Spec}
\DeclareMathOperator{\initial}{in}
\DeclareMathOperator{\Pic}{Pic}
\newtheorem{lemma}{Lemma}[section] 
\newtheorem{theorem}[lemma]{Theorem}
\newtheorem{corollary}[lemma]{Corollary}
\newtheorem{proposition}[lemma]{Proposition}
\newtheorem*{Thm4.1}{Theorem~\ref{t:sdbound}}
\newtheorem*{Thm4.10}{Theorem~\ref{t:hilbpolybound}}
\theoremstyle{definition}
\newtheorem{definition}[lemma]{Definition}
\newtheorem{algorithm}[lemma]{Algorithm}

\newtheorem{example}[lemma]{Example}
\newtheorem{remark}[lemma]{Remark}

\newtheorem*{acknowledgement}{Acknowledgements}
\theoremstyle{remark}
\newtheorem*{proof*}{Proof}
\renewcommand{\theequation}%
{\arabic{section}.\arabic{lemma}.\arabic{equation}}

\begin{document}

\title[Bounds on Multigraded Regularity]{Uniform bounds on multigraded
regularity}

\author{Diane Maclagan}
\address{Department of Mathematics \\ Stanford University \\ Stanford
  \\ CA 94305 \\ USA}
\email{maclagan@math.stanford.edu}

\author{Gregory G. Smith}
\address{Department of Mathematics \\ Barnard College \\
  Columbia University \\ New York \\ \newline NY 10027 \\ USA}
\email{ggsmith@math.columbia.edu}

\begin{abstract}
We give an effective uniform bound on the multigraded regularity of a
subscheme of a smooth projective toric variety $X$ with a given
multigraded Hilbert polynomial.  To establish this bound, we introduce
a new combinatorial tool, called a Stanley filtration, for studying
monomial ideals in the homogeneous coordinate ring of $X$.  As a
special case, we obtain a new proof of Gotzmann's regularity theorem.
We also discuss applications of this bound to the construction of
multigraded Hilbert schemes.
\end{abstract}

%\date{14 May 2003}

\maketitle

%%---------------------------------------------------------------
\section{Introduction}

Bounding the degree of the generators of a module or sheaf is a
central problem in commutative algebra and algebraic geometry.  The
modern approach to this problem concentrates on proving stronger
bounds involving Castelnuovo-Mumford regularity.  In fact,
Castelnuovo-Mumford regularity was introduced in \S14 of
\cite{Mumford} to bound the family of all projective subschemes having
a given Hilbert polynomial.  Following this appearance,
Castelnuovo-Mumford regularity has become a crucial ingredient in
bounding the degree of syzygies \cite{GLP} \cite{EL} and constructing
Hilbert schemes, Picard schemes and moduli spaces \cite{altmanKleiman}
\cite{Viehweg}.

The goal of this paper is to bound the multigraded Castelnuovo-Mumford
regularity (as defined in \cite{MaclaganSmith1}) of all subschemes of
a smooth projective toric variety $X$ that have a given multigraded
Hilbert polynomial.  To establish this bound, we work with saturated
monomial ideals in the homogeneous coordinate ring of $X$.  We
introduce a new combinatorial tool, called a Stanley filtration, for
studying monomial ideals.  Using an appropriate Stanley filtration, we
produce an effective bound for the multigraded regularity of an
individual ideal or family of ideals.  We also discuss applications of
our bound to the construction of multigraded Hilbert schemes.

Using ideals in the homogeneous coordinate ring $S$ to analyze
subschemes of $X$ has several advantages.  The $\mathbb{Z}^{r}$-graded
polynomial ring $S = \Bbbk[x_{1}, \dotsc, x_{n}]$, introduced in
\cite{Cox}, is intrinsic to the variety $X$.  By focusing on $X$
rather than a projective embedding of $X$, we reduce both the number
of variables and the total degree of the polynomials needed to
describe a subscheme.  When $\Pic(X) \neq \mathbb{Z}$, the
multigrading allows for stronger bounds on the equations defining a
subscheme.  Multigradings also produce a finer stratification of
subschemes of $X$.

The novel approach required for multigraded polynomial rings leads to
new insights in the standard graded case.  Indeed, when $X =
\mathbb{P}^{d}$, we obtain a new proof of Gotzmann's optimal bound on
the regularity of all subschemes having a given Hilbert polynomial.
Gotzmann's original proof \cite{gotzmann} relies on Macaulay's
characterization of the Hilbert function of an ideal in a standard
graded polynomial ring.  Since there is no version of Macaulay's
theorem for nonstandard gradings, the methods used in \cite{gotzmann}
do not apply in our situation.  In fact, there is typically no
lexicographic ideal in the homogeneous coordinate ring of $X$ (see
\cite{ACdN}) so we cannot expect a direct analogy of Macaulay's
result.  The alternative proof of Gotzmann's result given in
\cite{GreenGot}, also see \cite{Green} and \S4.3 in \cite{BH}, uses an
induction on a general hyperplane section.  Because a general
hypersurface is rarely a toric variety, this approach does not extend
to toric varieties.

The main combinatorial tool used in this paper is based on a Stanley
decomposition.  Given a monomial ideal $I$ in $S$, a Stanley
decomposition for $S/I$ is a set $\mathfrak{S}$ of pairs
$(\bm{x}^{\bm{u}}, \sigma)$ such that $S/I \cong
\bigoplus\nolimits_{(\bm{x}^{\bm{u}}, \sigma) \in \mathfrak{S}}
S_{\sigma} \bigl(-\deg( \bm{x}^{\bm{u}}) \bigr)$, where
$\bm{x}^{\bm{u}}$ is a monomial in $S$, $\sigma \subseteq \{ 1,
\dotsc, n\}$ and $S_{\sigma} = \Bbbk[x_{i} : i \in \sigma]$.  In other
words, if we identify the pair $(\bm{x}^{\bm{u}}, \sigma)$ with the
set $\{ \bm{x}^{\bm{u} + \bm{v}} \in S : \bm{x}^{\bm{v}} \in
S_{\sigma} \}$, then each monomial of $S$ not in $I$ belongs to a
unique pair $(\bm{x}^{\bm{u}}, \sigma)$.  It follows that a Stanley
decomposition expresses the multigraded Hilbert polynomial of $S/I$ as
a sum of the Hilbert polynomials for $S_{\sigma}$:
\begin{equation} \label{decomposition}
P_{S/I}(\bm{t})= \sum_{(\bm{x}^{\bm{u}}, \sigma) \in \mathfrak{S}}
P_{S_{\sigma}} \bigl( \bm{t} - \deg(\bm{x}^{\bm{u}}) \bigr) \, .
\end{equation}

\begin{example} \label{introexample}
Let $S = \Bbbk[x_{1}, \dotsc, x_{4}]$ have standard grading defined by
$\deg(x_{i}) = 1$ for $1 \leq i \leq 4$.  If $I = \langle
x_{1}x_{4}^{2}, x_{2}x_{4}^{2}, x_{3}x_{4}^{2} \rangle$ is an ideal in
$S$ then
\begin{align*}
&\text{\small $\bigl\{ (1, \{ 1, 2, 3 \} ), ( x_{4}, \{ 1, 2, 3\} ),
(x_{4}^{2}, \{ 4 \}) \bigr\}$} \text{ and } \\
&\text{\small $\bigr\{ (1, \{ 4 \}), (x_{3}, \{ 3 \} ), ( x_{3}x_{4},
\{ 4 \} ), (x_{2}, \{ 2, 3 \} ), (x_{2}x_{4}, \{ 2, 3 \}), (x_{1}, \{
1, 2, 3 \}), ( x_{1}x_{4}, \{ 1, 2, 3 \}) \bigr\}$}
\end{align*} 
are both Stanley decompositions for $S/I$.  Since the Hilbert
polynomial $P_{S_{\sigma}}(t)$ is simply the binomial coefficient
$\binom{t + |\sigma| -1}{|\sigma| +1}$, these Stanley decompositions
yield
\begin{align*}
P_{S/I}(t) &= \textstyle\binom{t+2}{2} + \binom{t+1}{2} +
\binom{t-2}{0} \\
&= \textstyle\binom{t}{0} + \binom{t-1}{0} + \binom{t-2}{0} +
\binom{t}{1} + \binom{t-1}{1} + \binom{t+1}{2} + \binom{t}{2} \, .
\end{align*}
\end{example}

We focus on a particular class of Stanley decompositions called
Stanley filtrations.  By definition, these are ordered sets $\{
(\bm{x}^{\bm{u}_{i}}, \sigma_{i} ) : 1 \leq i \leq m \}$ such that the
modules $M_{j} = S / \bigl(I + \langle \bm{x}^{\bm{u}_{j+1}}, \dotsc,
\bm{x}^{\bm{u}_{m}} \rangle \bigr)$ form a filtration $\Bbbk = M_{0}
\subset M_{1} \subset \dotsb \subset M_{m} = S/I$ with $M_{i}/M_{i-1}
= S_{\sigma_{i}}$.  The decompositions of Example~\ref{introexample}
are Stanley filtrations in the order presented.  We provide an
algorithm for finding Stanley filtrations.

Our first major theorem uses a Stanley filtration to give an effective
bound on the multigraded regularity.  Recall that bounding the
multigraded regularity of a module $M$ is equivalent to giving a
subset of $\reg(M) = \bigl\{ \bm{k} \in \mathbb{Z}^{r} : \text{$M$ is
$\bm{k}$-regular} \bigr\}$.  For more information on multigraded
regularity, we refer to \cite{MaclaganSmith1}.  Remarkably, our major
theorems use only the behavior of multigraded regularity in short
exact sequences and hence are independent of the precise definition of
multigraded regularity.

\begin{Thm4.1} 
Let $I$ be a monomial ideal in $S$.  If $\bigl\{ (
\bm{x}^{\bm{u}_{i}}, \sigma_{i}) : 1 \leq i \leq m \bigr\}$ is a
Stanley filtration for $S/I$, then $\bigcap_{i = 1}^{m} \bigl( \deg(
\bm{x}^{\bm{u}_{i}}) + \reg(S_{\sigma_{i}}) \bigr) \subseteq \reg
\left( S/I \right)$.
\end{Thm4.1}

By relating the sets $\sigma_{i}$ to the fan $\Delta$ defining $X$, we
can eliminate certain pairs from this intersection.

\begin{example}
Since $\reg(S) = \mathbb{N}$ for any standard graded polynomial ring,
the first Stanley filtration in Example~\ref{introexample} implies
that 
\[
\max\{ \deg(1), \deg(x_{4}), \deg(x_{4}^{2}) \} + \mathbb{N} \subseteq
\reg(S/I) \, .
\]  
The minimal free resolution of $S/I$ shows that this bound is sharp:
$\reg(S/I) = \bm{2} + \mathbb{N}$.
\end{example}

To study all subschemes of $X$ with a given multigraded Hilbert
polynomial, we use the combinatorial structure of $\Delta$ to focus on
a finite set of Stanley filtrations.  We also concentrate on ideals
that are saturated with respect to the irrelevant ideal $B$; see
Section~2.  Given a polynomial $P(\bm{t})$, we are most interested in
expressions of the form \eqref{decomposition} arising from our finite
set of Stanley filtrations.  This leads to an algorithm for finding
all $B$-saturated monomial ideals with multigraded Hilbert polynomial
$P(\bm{t})$.  We call the maximum number of summands in such an
expression for $P(\bm{t})$ the Gotzmann number.

To state our second major result, let $\widehat{\sigma}$ denote the
complement of $\sigma$ in $\{ 1, \dotsc, n\}$.  Identifying $\Pic(X)$
with $\mathbb{Z}^{r}$, we write $\scr{K} \subset \mathbb{Z}^{r}$ for
the semigroup of nef line bundles on $X$; see Section~2 for a
combinatorial description of $\scr{K}$.

\begin{Thm4.10} 
Let $I$ be any $B$-saturated ideal in $S$ and let $\bm{c} \in
\bigcap_{i=1}^{n} \bigl( \deg(x_{i}) + \scr{K} \bigr)$.  If $m$ is the
Gotzmann number for $P_{S/I}(\bm{t})$ then
$\bigcap\limits_{\widehat{\sigma} \in \Delta} \bigl( (m-1) \bm{c} +
\reg(S_{\sigma}) \bigr) \subseteq \reg(S/I)$.
\end{Thm4.10}

This theorem implies that for any $\bm{k} \in
\bigcap\nolimits_{\widehat{\sigma} \in \Delta} \bigl( (m-1) \bm{c} +
\reg(S_{\sigma}) \bigr)$ every the subscheme of $X$ having multigraded
Hilbert polynomial $P_{S/I}(\bm{t})$ is cut out by equations of
multidegree $\bm{k}$.  Specializing to $X = \mathbb{P}^{d}$, we
recover Gotzmann's regularity theorem; see Theorem~\ref{t:gotzmann}.

The structure of the paper is as follows.  The next section
establishes our notation for toric varieties, recalls the definition
of multigraded regularity from \cite{MaclaganSmith1} and collects the
basic properties of multigraded Hilbert polynomials.  In
Section~\ref{s:sds}, we develop the theory of Stanley decompositions
and filtrations.  The proofs of our major theorems are in
Section~\ref{s:boundandalgo}.  This section also contains the
algorithm for finding all $B$-saturated monomial ideals with a given
multigraded Hilbert polynomial.  In Section~\ref{s:standardcase}, we
restrict to the case $X = \mathbb{P}^{d}$ and show that multigraded
techniques provide a simple new proof of Gotzmann's regularity
theorem.  Finally, Section~\ref{s:hilbschemes} discusses the effective
construction of multigraded Hilbert schemes.

\begin{acknowledgement}
We thank Ezra Miller for the references on cleanness.  We also thank
Kristina Crona and Edwin O'Shea for their helpful comments on a
preliminary version of this paper.  Both authors were partially
supported by the Mathematical Sciences Research Institute in Berkeley,
CA.
\end{acknowledgement}

%%---------------------------------------------------------------
\section{Castelnuovo-Mumford Regularity and Hilbert Polynomials} 
\label{s:defns}

This section relates multigraded Hilbert polynomials to multigraded
Castelnuovo-Mumford regularity (as defined in \cite{MaclaganSmith1}).
Let $X$ be a smooth projective toric variety over a field $\Bbbk$
determined by a fan $\Delta$ in $\mathbb{R}^{d}$.  By numbering the
rays (one-dimensional cones), we identify $\Delta$ with a simplicial
complex on $[n] := \{ 1, \dotsc, n \}$.  We write $\bm{b}_{1}, \dotsc,
\bm{b}_{n}$ for the unique minimal lattice vectors generating the rays
and we assume that $\bm{b}_{1}, \dotsc, \bm{b}_{n}$ span
$\mathbb{R}^{d}$.  Set $r := n - d$ and fix an $(r \times n)$-matrix
$A = [ \bm{a}_{1} \dotsb \bm{a}_{n} ]$ such that there is a short
exact sequence
\begin{equation} \label{dualitySES}
0 \longrightarrow \mathbb{Z}^{d} \xrightarrow{\; [ \bm{b}_{1} \dotsb
\, \bm{b}_{n} ]^{\textsf{T}} \;} \mathbb{Z}^{n} \xrightarrow{\; [
\bm{a}_{1} \dotsb \, \bm{a}_{n} ] \;} \mathbb{Z}^{r} \longrightarrow 0
\, .
\end{equation}
Because $A$ is the Gale dual of the $(d \times n)$-matrix $[
\bm{b}_{1} \dotsb \bm{b}_{n} ]$, it is uniquely determined up to
unimodular (determinant $\pm 1$) coordinate transformations of
$\mathbb{Z}^{r}$.  Since $X$ is smooth, the Picard group of $X$ is
isomorphic to $\mathbb{Z}^{r}$.  The homogeneous coordinate ring of
$X$, introduced in \cite{Cox}, is the polynomial ring $S =
\Bbbk[x_{1}, \dotsc, x_{n}]$ with the $\mathbb{Z}^{r}$-grading defined
by $\deg(x_{i}) = \bm{a}_{i} \in \mathbb{Z}^{r}$ for $1 \leq i \leq
n$.  The combinatorial structure of $\Delta$ is encoded in the
irrelevant ideal $B = \bigl\langle \prod_{i \not\in \sigma} x_{i} :
\sigma \in \Delta \bigr\rangle$.  With these definitions, \cite{Cox}
proves that the category of coherent $\mathscr{O}_{X}$-modules is
equivalent to the category of finitely generated
$\mathbb{Z}^{r}$-graded $S$-modules modulo $B$-torsion modules.

\begin{example} 
When $X = \mathbb{P}^{d}$, the short exact sequence \eqref{dualitySES}
is
\[
0 \longrightarrow \mathbb{Z}^{d} \longrightarrow \mathbb{Z}^{d+1}
\xrightarrow{\; [ 1 \; 1 \, \dotsb \, 1 ] \;} \mathbb{Z}^{1}
\longrightarrow 0 \, .
\]  
Since $\deg(x_{i}) = \bm{a}_{i} = 1$ for $1 \leq i \leq n = d+1$, the
homogeneous coordinate ring $S = \Bbbk[x_{1}, \dotsc, x_{n}]$ is
simply the standard graded polynomial ring.  The irrelevant ideal $B$
is the unique graded maximal ideal $\langle x_{1}, \dotsc, x_{n}
\rangle$.
\end{example}

\begin{example} \label{e:Hirzebruch} 
If $X$ is the Hirzebruch surface (or rational scroll)
$\mathbb{F}_{\ell} = \mathbb{P} \big( \mathscr{O}_{\mathbb{P}^{1}}
\oplus \mathscr{O}_{\mathbb{P}^{1}}(\ell) \big)$, then the short exact
sequence \eqref{dualitySES} is
\[
0 \longrightarrow \mathbb{Z}^{2} \xrightarrow{\; \left[
\begin{smallmatrix} 
1 & 0 & -1 & 0 \\ 
0 & 1 & \ell & -1 
\end{smallmatrix} 
\right]^{\textsf{T}} \; } \mathbb{Z}^{4} \xrightarrow{\; \left[   
\begin{smallmatrix} 
1 & - \ell & 1 & 0 \\ 
0 & 1 & 0 & 1   
\end{smallmatrix} \right] 
\;} \mathbb{Z}^{2} \longrightarrow 0 \, .
\]
Figure~\ref{f:Hirzfig} illustrates the fan and grading $A$ when $\ell =
2$.
\begin{figure}[ht]
\psfrag{1}{\scriptsize $1$}
\psfrag{2}{\scriptsize $2$}
\psfrag{3}{\scriptsize $3$}
\psfrag{4}{\scriptsize $4$}
\psfrag{A}{$A$}
\psfrag{B}{$\Delta$}
\epsfig{file=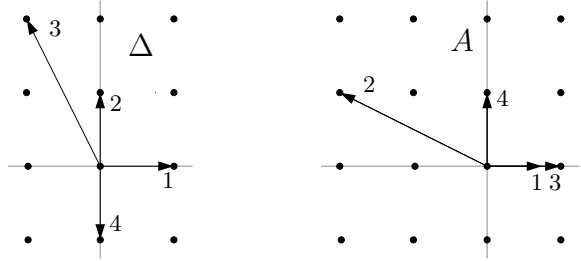, width=3in}
\caption{The fan and grading for $\mathbb{F}_{2}$. \label{f:Hirzfig}}
\end{figure}
The homogeneous coordinate ring $S = \Bbbk[x_{1}, x_{2}, x_{3},
x_{4}]$ has the $\mathbb{Z}^{2}$-grading induced by $\deg(x_{1}) =
\left[
\begin{smallmatrix} 1 \\ 0 \end{smallmatrix} \right]$,
$\deg(x_{2}) = 
\left[ \begin{smallmatrix} - \ell \\ 1 \end{smallmatrix} \right]$,
$\deg(x_{3}) = \left[
\begin{smallmatrix} 1 \\ 0 \end{smallmatrix} \right]$,
$\deg(x_{4}) = \left[ 
\begin{smallmatrix} 0 \\ 1 \end{smallmatrix} 
\right]$ and $B = \langle x_{1}x_{2}, x_{2}x_{3}, x_{3}x_{4},
x_{1}x_{4} \rangle = \langle x_{1}, x_{3} \rangle \cap \langle x_{2},
x_{4} \rangle$.
\end{example}

The combinatorial structure of $\Delta$ also gives rise to an
important subsemigroup of $\mathbb{Z}^{r}$.  We write $\mathbb{N}
A_{\sigma} := \{ \sum_{i \in \sigma} \lambda_{i} \bm{a}_{i} :
\lambda_{i} \in \mathbb{N} \}$ for the affine semigroup generated by
the set $\{ \bm{a}_{i} : i \in \sigma \}$.  For $\sigma \subseteq
[n]$, let $\widehat{\sigma}$ denote the complement of $\sigma$ in
$[n]$.  The semigroup $\scr{K}$ is $\bigcap_{\sigma \in \Delta}
\mathbb{N} A_{\widehat{\sigma}}$.  Since $X$ is projective, $\scr{K}$
is the set of integral points of an $r$-dimensional pointed cone in
$\mathbb{R}^{r}$.  Geometrically, elements in $\scr{K}$ correspond to
numerically effective (nef) line bundles on $X$.  As \cite{cox2}
indicates, $\scr{K} \otimes_{\mathbb{Z}} \mathbb{R}$ is the closure of
the K\"{a}hler cone of $X$.  The dual of the K\"{a}hler cone is the
Mori cone of effective $1$-cycles modulo numerical equivalence.

\begin{example} \label{e:examplecones}
When $X = \mathbb{P}^{d}$, the semigroup $\scr{K} = \mathbb{N}$.  If
$X = \mathbb{F}_{\ell}$ (with $A$ chosen as in
Example~\ref{e:Hirzebruch}), then $\scr{K} = \mathbb{N}^{2}$.  In
general, the structure of $\scr{K}$ can be much more complicated; see
Example~2.8 in \cite{MaclaganSmith1}.
\end{example}

The next result illustrates the connection between the irrelevant
ideal $B$ and the semigroup $\scr{K}$.  For $\sigma \subseteq [n]$,
let $P_{\sigma}$ be the prime ideal $\langle x_{i} : i \not\in \sigma
\rangle$ and let $S_{\sigma}$ be the (smaller) polynomial ring
$\Bbbk[x_{i} : i \in \sigma] \cong S/P_{\sigma}$.

\begin{lemma} \label{l:Bsaturated}
A monomial ideal $I$ in $S$ is $B$-saturated if and only if every
associated prime $P_{\sigma}$ of $I$ satisfies $\widehat{\sigma} \in
\Delta$ \text{$($or equivalently $\scr{K} \subseteq \mathbb{N}
A_{\sigma})$}.
\end{lemma}

\begin{proof}
Let $I = \bigcap_{\sigma} Q_{\sigma}$ be an irredundant primary
decomposition for $I$ where the ideal $Q_{\sigma}$ is
$P_{\sigma}$-primary.  It follows that $(I : B^{\infty}) =
\bigcap_{\sigma} (Q_{\sigma} : B^{\infty})= \bigcap_{\sigma}
\bigcap_{\tau \in \Delta} \bigl(Q_{\sigma} : (\prod_{i \not \in \tau}
x_i)^{\infty} \bigr)$.  Now $\bigl( Q_{\sigma} : (\prod_{i \not \in
\tau} x_i)^{\infty} \bigr)$ equals $S$ if $\widehat{\sigma} \cap
\widehat{\tau} \neq \emptyset$ and equals $Q_{\sigma}$ otherwise.
Since $\widehat{\sigma} \cap \widehat{\tau} = \emptyset$ is
equivalent to $\widehat{\sigma} \subseteq \tau$, we have
$(I:B^{\infty})= \bigcap_{\widehat{\sigma} \in \Delta} Q_{\sigma}$.
Therefore, $(I:B^{\infty})=I$ if and only if every associated prime
$P_{\sigma}$ satisfies $\widehat{\sigma} \in \Delta$.  The equivalent
condition follows immediately from the definition of $\scr{K}$.
\end{proof}

Throughout this paper, $M$ denotes a finitely generated
$\mathbb{Z}^{r}$-graded $S$-module.  We refer to \cite{BrodmannSharp}
for background information on local cohomology.  The module $M$ is
$B$-torsion if $M = H_{B}^{0}(M) = \bigcup_{j \in \mathbb{N}} (0 :_{M}
B^{j})$.  At the other extreme, $M$ is $B$-torsion-free if
$H_{B}^{0}(M) = 0$.  For an ideal $I \subseteq S$, the module $S/I$ is
$B$-torsion-free if and only if $I$ is $B$-saturated; $(I:B^{\infty})
= I$.

\begin{remark} \label{r:Btorsion}
The semigroup $\scr{K}$ also has a useful algebraic interpretation.
If $\Bbbk$ is infinite and $M$ is $B$-torsion-free, then
Proposition~3.1 in \cite{MaclaganSmith1} shows that for any $\bm{k}
\in \scr{K}$ there is a nonzerodivisor $f \in S_{\bm{k}}$ on $M$.
\end{remark}

Let $\scr{C} = \{ \bm{c}_{1}, \dotsc, \bm{c}_{e} \}$ be the unique
minimal Hilbert basis of $\scr{K}$.  By definition, $\scr{C}$ is the
minimal subset of $\scr{K}$ such that every element in $\scr{K}$ is a
nonnegative integral combination of the $\bm{c}_{j}$; see
\S~\!\!IV.16.4 in \cite{schrijver}.  We recall the definition of
multigraded Castelnuovo-Mumford regularity introduced in
\cite{MaclaganSmith1}.

\begin{definition} \label{d:localreg}
For $\bm{k} \in \mathbb{Z}^{r}$, the module $M$ is
\emph{$\bm{k}$-regular} if the following conditions are satisfied:
\begin{enumerate}
\item $H_{B}^{i}(M)_{\bm{p}} = 0$ for all $i \geq 1$ and all $\bm{p}
\in \bigcup ( \bm{k} - \lambda_{1} \bm{c}_{1} - \dotsb - \lambda_{e}
\bm{c}_{e} + \scr{K})$ where the union is over all $\lambda_{1},
\dotsc, \lambda_{e} \in \mathbb{N}$ such that $\lambda_{1} + \dotsb +
\lambda_{e} = i-1$;
\item $H_{B}^{0}(M)_{\bm{p}} = 0$ for all $\bm{p} \in \bigcup_{1 \leq
j \leq e} ( \bm{k} + \bm{c}_{j} + \scr{K})$.
\end{enumerate}
The \emph{regularity of $M$}, denoted by $\reg(M)$, is the set $\{
\bm{k} \in \mathbb{Z}^{r} : \text{$M$ is $\bm{k}$-regular} \}$.
\end{definition}

In this paper, we exploit two properties of multigraded
Castelnuovo-Mumford regularity.  Firstly, if $I$ is an ideal in $S$
and $\bm{k} \in \reg(I)$ then the subscheme of $X$ defined by $I$ is
cut out by equations of multidegree $\bm{k}$; see Theorem~6.9 in
\cite{MaclaganSmith1}.  The second property allows us to focus on
monomial ideals by relating the regularity of an ideal with its
initial ideal.  We write $\initial(I)$ for the initial ideal of $I$
with respect to some monomial order.  The following proposition, which
is well-known for $X = \mathbb{P}^{d}$, appears as Proposition~6.13 in
\cite{MaclaganSmith1}.

\begin{proposition} \label{p:usc}
If $I$ is an ideal in $S$, then $\reg \bigl( S/\initial(I) \bigr)
\subseteq \reg(S/I)$.  Moreover, if $I$ is $B$-saturated and $J =
\bigl( \initial(I) : B^{\infty} \bigr)$ then $\reg(S/J) \subseteq
\reg(S/I)$. \qed
\end{proposition}

Next, we turn our attention to multigraded Hilbert polynomials.  The
multigraded Hilbert function $H(M, \bm{t})$ equals the dimension of
the degree $\bm{t}$ homogeneous component of a
$\mathbb{Z}^{r}$\nobreakdash-graded module $M$.  As in the standard
graded case, $H(M, \bm{t})$ ``eventually'' agrees with a polynomial.
To prove this, we first consider the multigraded Hilbert function of
the ring $S$.

\begin{lemma} \label{l:polynomial}
If $S$ is the coordinate ring of a smooth toric variety then the
Hilbert function $H(S, \bm{t}) = \dim_{\Bbbk} S_{\bm{t}}$ agrees with
a polynomial for all $\bm{t} \in \scr{K}$.
\end{lemma}

\begin{proof}
Since the monomials of degree $\bm{t}$ form a basis for the
$\Bbbk$-vector space $S_{\bm{t}}$, the Hilbert function $H(S, \bm{t})$
is a vector partition function.  This means $H(S, \bm{t})$ equals the
number of ways a vector $\bm{t} \in \mathbb{Z}^{r}$ can be written as
a sum of $\bm{a}_{1}, \dotsc, \bm{a}_{n}$.  The chamber complex of $\{
\bm{a}_{1}, \dotsc, \bm{a}_{n} \}$ is a polyhedral subdivision of
$\pos\{ \bm{a}_{1}, \dotsc, \bm{a}_{n} \}$.  It is defined to be the
common refinement of the simplicial cones $\pos\{ \bm{a}_{i} : i
\not\in \sigma \}$ where $\sigma \in \Delta$.  Hence, the cone
$\scr{K} \otimes_{\mathbb{Z}} \mathbb{R}$ is a chamber (maximal cell)
in the chamber complex.  From \cite{SturmfelsVPF}, we know that vector
partition functions are piecewise quasi-polynomials on the chamber
complex.  Therefore, $H(S, \bm{t})$ is a quasi-polynomial on
$\scr{K}$.

To complete the proof, we show that the period of this
quasi-polynomial is one.  We write $[ \bm{a}_{i} : i \not\in \sigma ]$
for the submatrix of $A$ consisting of those columns indexed by
$\sigma$.  From \cite{SturmfelsVPF}, we know that the period of the
quasi-polynomial is at most the least common multiple of $\det [
\bm{a}_{i} : i \not\in \sigma ]$ where $\sigma$ is a facet in
$\Delta$.  By renumbering (if necessary) the $\bm{b}_{i}$, we may
assume that $\sigma = \{ 1, \dotsc, d\} \in \Delta$.  Recall that $X$
is smooth if and only if $\det [ \bm{b}_{i} : i \in \sigma] = \pm 1$
for all facets $\sigma \in \Delta$; see \S2.1 in \cite{fulton}.
Hence, there exists a unimodular change of coordinates such that
$\bm{b}_{i} = \bm{e}_{i}$ for all $i \in \sigma$ where $\bm{e}_{i}$ is
the $i$th standard basis vector.  In other words, $[ \bm{b}_{1} \dotsb
\bm{b}_{n} ]$ is the block matrix $[ I_{d} \, | \, V_{\sigma} ]$ where
$V_{\sigma}$ is a $(d \times r)$-matrix.  The Gale dual of this
configuration is $[ V_{\sigma}^{\textsf{T}} \, | \, -I_{r} ]$.
Because the Gale dual is determined up to unimodular transformation,
we have $\det [ \bm{a}_{i} : i \not \in \sigma] = \pm \det[ I_{r} ] =
\pm 1$ for all facets $\sigma \in \Delta$.
\end{proof}

Algorithms for computing $P_{S}(\bm{t})$ have been implemented in the
software package \texttt{LattE}; see \cite{LattE}.

\begin{example} \label{e:hilbertpoly}
When $X = \mathbb{P}^{d}$, we have $P_{S}(t) = \binom{t+d}{d}$.  If $X
= \mathbb{F}_{\ell}$, then we have $P_{S}(t_{1}, t_{2}) = t_{1}t_{2} +
\bigl( \tfrac{\ell}{2} \bigr) t_{2}^{2} + t_{1} + \bigl(
\tfrac{\ell+2}{2} \bigr) t_{2} + 1$.
\end{example}

Using Lemma~\ref{l:polynomial}, we show that the multigraded Hilbert
function of a module $M$ agrees with a polynomial for values of
$\bm{t}$ sufficiently far into the interior of $\scr{K}$.

\begin{proposition} \label{p:hilbpolyexists}
There exists a unique polynomial $P_{M}(\bm{t}) \in \mathbb{Q}[t_{1},
\dotsc, t_{r}]$ such that $P_{M}(\bm{t}) = H(M, \bm{t})$ for all
$\bm{t}$ in a finite intersection of translates of $\scr{K}$.  In
particular, $H(M,\bm{t})$ agrees with $P_{M}(\bm{t})$ for all $\bm{t}$
sufficiently far from the boundary of $\scr{K}$.
\end{proposition}

\begin{proof}
If $0 \longrightarrow \bigoplus_{j} S(- \bm{q}_{p, j}) \longrightarrow
\dotsb \longrightarrow \bigoplus_{j} S(- \bm{q}_{0, j})$ is the
minimal free resolution of $M$, then $H(M, \bm{t}) = \sum_{i = 0}^{p}
\sum_{j} (-1)^{i} H(S, \bm{t} - \bm{q}_{i,j})$.  It follows from
Lemma~\ref{l:polynomial} that $H(M, \bm{t})$ is a polynomial for all
$\bm{t} \in \bigcap_{i,j} ( \bm{q}_{i,j} + \scr{K})$.  Since
$\bm{q}_{i,j} \in \mathbb{Z}^{r}$ and $\scr{K}$ corresponds to the
lattice points in an $r$-dimensional cone in $\mathbb{R}^{r}$, this
intersection is nonempty.
\end{proof}

\begin{definition}
The polynomial $P_{M}(\bm{t})$ in Proposition~\ref{p:hilbpolyexists}
is called the \emph{multigraded Hilbert polynomial} of $M$.
\end{definition}

The multigraded Hilbert polynomial of a $B$-torsion module is
especially simple.

\begin{lemma} \label{l:torsionpoly}
If $M$ is a $B$-torsion module then $P_{M}(\bm{t}) = 0$.
\end{lemma}

\begin{proof}
Since $M$ is a finitely generated $B$-torsion module, there is a $j
\gg 0$ such that $B^{j} M = 0$.  We first prove that there exists a
$\bm{k} \in \scr{K}$ such that if $\bm{p} \in \bm{k} + \scr{K}$ then
every monomial in $S_{\bm{p}}$ belongs to $B^{j}$.  Choose an element
$\bm{c}$ which lies in the interior of $\scr{K}$.  If
$\bm{x}^{\bm{u}}$ is a monomial in $S$ and $\deg(\bm{x}^{\bm{u}}) \in
\bm{c} + \scr{K}$, then Lemma~2.4 in \cite{MaclaganSmith1} shows that
$\bm{x}^{\bm{u}} \in B$.  Suppose that $\bm{x}^{\bm{v}} \in S$ such
that $\deg(\bm{x}^{\bm{v}}) = n\bm{c} + \bm{c}'$ for $\bm{c}' \in
\scr{K}$.  Caratheodory's Theorem (Proposition~1.15 in \cite{Ziegler})
implies that $\bm{v} = \lambda_{1} \bm{u}_{1} + \dotsb + \lambda_{n}
\bm{u}_{n} + \bm{w}$ for some $\bm{u}_{1}, \dotsc, \bm{u}_{n} \in
\mathbb{N}^{n}$ satisfying $\deg( \bm{x}^{\bm{u}_{i}}) = \bm{c}$ for
$1 \leq i \leq n$, $\bm{w} \in \mathbb{R}_{\geq 0}^{n}$ satisfying $A
\bm{w} = \bm{c}'$ and some $\lambda_{1}, \dotsc, \lambda_{n} \in
\mathbb{R}_{\geq 0}$ satisfying $\lambda_{1} + \dotsb + \lambda_{n} =
n$.  It follows that there is an $i \in [n]$ such that $\lambda_{i}
\geq 1$ and hence $\bm{x}^{\bm{v}}$ is divisible by
$\bm{x}^{\bm{u}_{i}}$.  Therefore, if we set $\bm{k} := (j+n) \bm{c}$
then $\deg(\bm{x}^{\bm{v}}) \in \bm{k} + \scr{K}$ implies that
$\bm{x}^{\bm{v}} \in B^{j}$.

To complete the proof, we show that $M_{\bm{t}} = 0$ for all $\bm{t}$
sufficiently far into the interior of $\scr{K}$.  Let $f_{1}, \dotsc,
f_{h}$ be generators of $M$.  Our choice of $\bm{k}$ guarantees that
$M_{\bm{t}} = 0$ for all $\bm{t} \in \bigcap_{1 \leq i \leq h} \bigl(
\deg(f_{i}) + \bm{k} + \scr{K} \bigr)$.  Since elements in $\scr{K}$
are lattice points in a full-dimensional cone, the elements in this
intersection are the lattice points in a translation of the same cone.
We conclude that $P_{M}(\bm{t}) = 0$.
\end{proof}

More generally, the multigraded Hilbert polynomial of a module is
independent of $B$-torsion.

\begin{lemma} \label{l:saturationpoly}
If $\overline{M} := M / H_{B}^{0}(M)$ then $P_{M}(\bm{t}) =
P_{\overline{M}}(\bm{t})$.  In particular, if $I \subseteq S$ is an
ideal then $S/I$ and $S / (I : B^{\infty})$ have the same Hilbert
polynomial.
\end{lemma}

\begin{proof}
Since $H_{B}^{0}(M)$ is a $B$-torsion module,
Lemma~\ref{l:torsionpoly} shows that its multigraded Hilbert
polynomial equals $0$.  Hence, the short exact sequence 
\[
0 \longrightarrow H_{B}^{0}(M) \longrightarrow M \longrightarrow
\overline{M} \longrightarrow 0
\] 
implies that $P_{M}(\bm{t}) = P_{\overline{M}}(\bm{t})$.  Because
$H_{B}^{0}(S/I) = (I : B^{\infty})/I$, the second assertion is a
special case of the first part.
\end{proof}

The following result connects Hilbert functions with local cohomology
modules.  The special case in which $S$ has the standard grading can
be found in \S4.4 of \cite{BH}.

\begin{proposition} \label{p:funcpolydifference}
We have
\begin{equation} \label{funcpolydifference:1}
H(M, \bm{t}) - P_{M}(\bm{t}) = \sum_{i = 0}^{d} (-1)^{i} \dim_{\Bbbk}
H_{B}^{i}(M)_{\bm{t}} \qquad \text{for all $\bm{t} \in
\mathbb{Z}^{r}$.}
\end{equation}
\end{proposition}

\begin{proof}
We proceed by induction on $\dim M$.  If $\dim M = 0$, then $M$ is
artinian.  Hence, $M$ is a $B$-torsion module and we have $M =
H_{B}^{0}(M)$ and $H_{B}^{i}(M) = 0$ for all $i \geq 1$.  Since
Lemma~\ref{l:torsionpoly} shows that $P_{M}(\bm{t}) = 0$, the
assertion follows.

Assume $\dim M > 0$.  Since both sides of \eqref{funcpolydifference:1}
change by $\dim_{\Bbbk} H_{B}^{0}(M)_{\bm{t}}$ when $M$ is replaced by
$M / H_{B}^{0}(M)$, we may assume that $M$ is $B$-torsion-free.
Because extension of the base field commutes with the formation of
local cohomology, we may also assume that $\Bbbk$ is infinite.  Choose
$\bm{k} \in \scr{K}$.  Remark~\ref{r:Btorsion} implies there is a
nonzerodivisor $f$ on $M$ with $f \in S_{\bm{k}}$.  Hence, $\dim M/f M
< \dim M$ and there is a short exact sequence
\begin{equation} \label{nzdseq}
0 \longrightarrow M(- \bm{k}) \xrightarrow{\;\; f \;\;} M
\longrightarrow M/f M \longrightarrow 0 \, .
\end{equation}  
Set $H_{M}'(\bm{z}) = \sum\nolimits_{\bm{t} \in \mathbb{Z}^{r}} \big(
H(M, \bm{t}) - P_{M}(\bm{t}) \big) \, \bm{z}^{\bm{t}}$ and
\[
H_{M}''(\bm{z}) = \sum_{\bm{t} \in \mathbb{Z}^{r}} \left( \sum_{i =
0}^{d} (-1)^{i} \dim_{\Bbbk} H_{B}^{i}(M)_{\bm{t}} \right)
\bm{z}^{\bm{t}} \, .
\]
With this notation, it suffices to prove that $H_{M}'(\bm{z}) =
H_{M}''(\bm{z})$.  From \eqref{nzdseq}, it follows that $H(M /
fM,\bm{t}) = H(M, \bm{t}) - H(M, \bm{t} - \bm{k})$.  Combining this
with Proposition~\ref{p:hilbpolyexists}, we deduce that $P_{M /
fM}(\bm{t}) = P_{M}(\bm{t}) - P_{M}(\bm{t}-\bm{k})$ for all $\bm{t}$
sufficiently far into the interior of $\scr{K}$ and thus for all
$\bm{t}$.  Hence, we have $H_{M / f M}'(\bm{z}) = (1 -
\bm{z}^{\bm{k}}) H_{M}'(\bm{z})$.  On the other hand, the long exact
sequence associated to \eqref{nzdseq} shows that 
\[
\sum_{i = 0}^{d} (-1)^{i} \dim_{\Bbbk} H_{B}^{i}(M /fM)_{\bm{t}} =
\sum_{i = 0}^{d} (-1)^{i} \bigl( \dim_{\Bbbk} H_{B}^{i}(M)_{\bm{t}} -
\dim_{\Bbbk} H_{B}^{i}(M)_{\bm{t} - \bm{k}} \bigr) \, .
\]
Therefore, we have $H_{M/fM}''(\bm{z}) = (1 - \bm{z}^{\bm{k}})
H_{M}''(\bm{z})$.  Since the induction hypothesis yields $H_{M / f
M}'(\bm{z}) = H_{M/fM}''(\bm{z})$, we conclude that $H_{M}'(\bm{z}) =
H_{M}''(\bm{z})$.
\end{proof}

\begin{corollary} \label{c:regmeanspoly}
If $M$ is $\bm{k}$-regular then the Hilbert function $H(M, \bm{t})$
agrees with the Hilbert polynomial $P_{M}(\bm{t})$ for all values
$\bm{t} \in \bm{k} + \scr{K}$ with $\bm{t} \neq \bm{k}$.
\end{corollary}

\begin{proof}
If $M$ is $\bm{k}$-regular, then $H_{B}^{i}(M)_{\bm{t}} =
0$ for all $i \geq 0$ and all $\bm{t} \in \bm{k} +
\scr{K}$ with $\bm{t} \neq \bm{k}$.  Hence, the claim
follows from Proposition~\ref{p:funcpolydifference}.
\end{proof}

Multigraded Hilbert polynomials also have a geometric description,
which is attributed to Snapper in \cite{Kleiman}.  Let
$\mathscr{O}_{X}(\bm{t})$ be the line bundle on $X$ corresponding to
$\bm{t} \in \mathbb{Z}^{r}$ and let $\mathscr{F}$ be the
$\mathscr{O}_{X}$\nobreakdash-module associated to the $S$-module $M$.
Since equation (6.3.1) in \cite{MaclaganSmith1} indicates that
\[
\dim_{\Bbbk} H^{0} \bigl( X, \mathscr{F} \otimes
\mathscr{O}_{X}(\bm{t}) \bigr) = H(M, \bm{t}) - \dim_{\Bbbk}
H_{B}^{0}(M)_{\bm{t}} + \dim_{\Bbbk} H_{B}^{1}(M)_{\bm{t}} \, ,
\]
Proposition~\ref{p:funcpolydifference} implies that
\[
P_{M}(\bm{t}) = \chi \bigl( \mathscr{F} \otimes
\mathscr{O}_{X}(\bm{t}) \bigr) = \sum\limits_{i =1}^{d} (-1)^{i} \dim
H^{i} \bigl( X, \mathscr{F} \otimes \mathscr{O}_{X}( \bm{t}) \bigr) \,
.
\]

For a finite set of points, the connection between multigraded
regularity and multigraded Hilbert polynomials is particularly
elegant.

\begin{example}
Let $I$ be the $B$-saturated ideal corresponding to a finite set of
points on $X$.  Proposition~6.7 in \cite{MaclaganSmith1} shows that
$\reg(S/I)$ is exactly the subset of $\mathbb{Z}^{r}$ for which the
Hilbert function $H(S/I, \bm{t})$ equals the Hilbert polynomial
$P_{S/I}(\bm{t})$.
\end{example}

%%---------------------------------------------------------------
\section{Stanley Decompositions and Filtrations} \label{s:sds}

In this section we introduce the key combinatorial tool used in this
paper.  We restrict our focus to a monomial ideal $I$ in the
polynomial ring $S$, and introduce the notion of a Stanley
decomposition for $S/I$.  This is a partition of the monomials of $S$
not in $I$ into sets each of which corresponds to the monomials in a
smaller polynomial ring.

\begin{definition} \label{d:sddefn}
If $\bm{x}^{\bm{u}} \in S$ and $\sigma \subseteq [n]$, the pair $(
\bm{x}^{\bm{u}}, \sigma)$ denotes the set of all monomials in $S$ of
the form $\bm{x}^{\bm{v} + \bm{u}}$ where $\supp(\bm{v}) := \{ i :
v_{i} \neq 0 \} \subseteq \sigma$.  A \emph{Stanley decomposition} for
$S/I$ is a set $\mathfrak{S}$ of pairs $(\bm{x}^{\bm{u}}, \sigma)$
such that
\[
S/I \cong \bigoplus_{(\bm{x}^{\bm{u}}, \sigma) \in
\mathfrak{S}} S_{\sigma}(-A \bm{u}) \, ,
\]
where $S_{\sigma} := \Bbbk[x_{i} : i \in \sigma]$.  In other words,
each monomial of $S$ not in $I$ belongs to a unique pair
$(\bm{x}^{\bm{u}}, \sigma)$ in the Stanley decomposition.
\end{definition}

A Stanley decomposition $\mathfrak{S}$ for $S/I$ also gives a primary
decomposition of $I$:
\[
I = \bigcap_{(\bm{x}^{\bm{u}}, \sigma) \in \mathfrak{S}} \langle
x_{i}^{u_{i} + 1} : i \not\in \sigma \rangle \, .
\]
This is typically not the unique irreducible irredundant primary
decomposition of $I$.  Stanley decompositions are inspired by
\cite{StanleySD} and algorithmically defined in \cite{SturmfelsWhite};
also see \cite{HoTh1}, \cite{apel}.  Both \cite{StanleySD} and
\cite{apel} require the extra condition that $|\sigma|$ should be at
least the depth of $I$.

\begin{example} \label{e:sdegs}
If $I = \langle x_{1}^{2}x_{2}, x_{1}x_{2}^{2} \rangle \subset S =
\Bbbk[x_{1}, x_{2}]$, then
\begin{enumerate}
\item $\bigl\{ \bigl( 1, \{ 1 \} \bigr), \bigl( x_{2},\{ 2 \} \bigr),
\bigl( x_{1}x_{2}, \emptyset \bigr) \bigr\}$,
\item$\bigl\{ \bigl( 1, \{ 2 \} \bigr), \bigl( x_{1}, \{ 1 \} \bigr),
\bigl( x_{1}x_{2}, \emptyset \bigr) \bigr\}$ and
\item $\bigl\{ \bigl( 1, \emptyset \bigr), \bigl( x_{1}, \{ 1 \}
\bigr), \bigl( x_{2}, \emptyset \bigr), \bigl( x_{2}^{2}, \{ 2 \}
\bigr), \bigl( x_{1}x_{2}, \emptyset \bigr) \bigr\}$ 
\end{enumerate}
are three distinct Stanley decomposition for $S/I$.  These are
illustrated in Figure~\ref{f:firstsdpic}.
\begin{figure}[ht] 
\psfrag{x1}{\scriptsize $x_1$}
\psfrag{x2}{\scriptsize $x_2$}
\psfrag{1}{\textbf{1}}
\psfrag{2}{\textbf{2}}
\psfrag{3}{\textbf{3}}
\epsfig{file=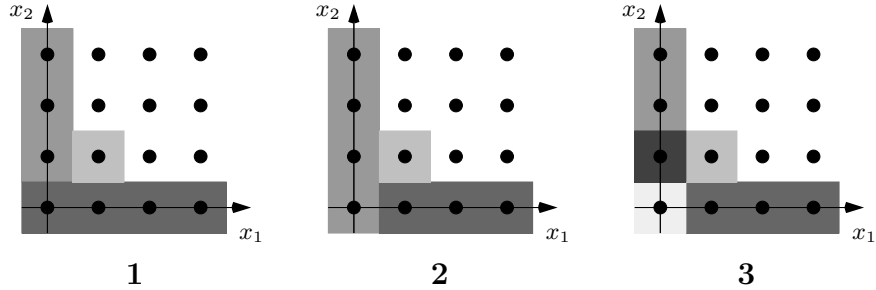, width=4.5in}
\caption{Stanley decompositions for $\langle x_1^2x_2,x_1x_2^2
\rangle$. \label{f:firstsdpic}}
\end{figure}
\end{example}

Stanley decompositions are closely related to standard pairs.  See
\cite{STV} for the origin of the notation $(\bm{x}^{\bm{u}},\sigma)$
and more details.  Standard pairs enjoy the following property: if
$(\bm{x}^{\bm{u}}, \sigma)$ is a standard pair of $I$ then $P_{\sigma}
= \langle x_{i} : i \not\in \sigma \rangle$ is an associated prime of
$I$.  In contrast, not all ideals have a Stanley decomposition where
every $\sigma$ corresponds to an associated prime.  

\begin{example} 
If $I = \langle x_{1}, x_{2} \rangle \cap \langle x_{3}, x_{4} \rangle
= \langle x_{1}x_{3}, x_{1}x_{4}, x_{2}x_{3}, x_{2}x_{4} \rangle$ is a
monomial ideal in the ring $S = \Bbbk[x_{1}, x_{2}, x_{3}, x_{4}]$,
then $\bigl\{ \bigl( 1, \{ 1, 2 \} \bigr), \bigl( x_{4}, \{ 4 \}
\bigr), \bigl( x_{3}, \{ 3, 4 \} \bigr) \bigr\}$ is a Stanley
decomposition for $S/I$ where $\{ 4 \}$ does not correspond to an
associated prime of $I$.  One easily verifies for this ideal that
every Stanley decomposition for $S/I$ has a pair $(\bm{x}^{\bm{u}},
\sigma)$ for which $\sigma$ does not correspond to an associated prime.
\end{example}

The paper \cite{clean} studies the special case when $S/I$ has a
Stanley decomposition in which each $\sigma$ corresponds to a minimal
associated prime of $I$.  Decompositions with this property are called
\emph{clean}.

One way to construct a Stanley decomposition is to make repeated use
of the short exact sequence
\begin{equation} \label{ses}
0 \longrightarrow S(- \bm{a}_{i})/(I:x_{i}) \xrightarrow{\;\; x_{i}
\;\;} S/I \longrightarrow S/(I + \langle x_{i} \rangle)
\longrightarrow 0 \, ,
\end{equation}
where $x_{i}$ is any variable.  More explicitly, we have following
algorithm.  A special case of this algorithm is implicit in the proof
of Lemma~2.4 in \cite{SturmfelsWhite}.

\begin{algorithm} \label{a:comma-colon}
Given a monomial ideal $I$ in the polynomial ring $S$ with $I \neq S$,
the following algorithm computes a Stanley decomposition for $S/I$.
\begin{enumerate}
\item (Base case) If $I$ is a prime ideal, then let $\sigma$
correspond to the set of variables not in $I$ and output $\bigl\{ (1,
\sigma) \bigr\}$.
\item (Choose variable) If $I$ is not prime then choose a variable
$x_{\ell} \in S$ that is a proper divisor of a minimal generator of
$I$.
\item (Recursion) Compute a Stanley decomposition $\bigl\{ (
\bm{x}^{\bm{u}}, \tau) \bigr\}$ for $S/ ( I + \langle
x_{\ell} \rangle)$ and a Stanley decomposition $\bigl\{ (
\bm{x}^{\bm{v}}, \sigma) \bigr\}$ for $S/(I : x_{\ell})$.
Output $\bigl\{ (\bm{x}^{\bm{u}}, \tau) \bigr\} \cup
\bigl\{ (\bm{x}^{\bm{v}} x_{\ell}, \sigma) \bigr\}$.
\end{enumerate}
\end{algorithm}

\begin{proof}[Proof of Correctness]
A monomial ideal is prime if and only if it is generated by a subset
of the variables.  Hence, when $I$ is prime and $\sigma$ corresponds
to the set of variables not in $I$, the set $\bigl\{ (1, \sigma)
\bigr\}$ is a Stanley decomposition for $S/I$.  On the other hand, if
$I$ is not prime, then there exists a variable $x_{\ell}$ that is a
proper divisor of a minimal generator of $I$.  From \eqref{ses}, we
see that a monomial not in $I$ corresponds to either a monomial not
in $S/(I + \langle x_{\ell} \rangle)$ or $x_{\ell}$ times a monomial
not in $S/(I:x_{\ell})$.  Thus, if $\bigl\{ (\bm{x}^{\bm{u}}, \tau)
\bigr\}$ is a Stanley decomposition for $S/(I + \langle x_{\ell}
\rangle)$ and $\bigl\{ (\bm{x}^{\bm{v}}, \sigma) \bigr\}$ is a Stanley
decomposition for $S/(I : x_i)$, then $\bigl\{ (\bm{x}^{\bm{u}}, \tau)
\bigr\} \cup \bigl\{ (\bm{x}^{\bm{v}}x_{\ell}, \sigma) \bigr\}$ is a
Stanley decomposition for $S/I$.  Finally, the algorithm terminates
because $S$ is a noetherian ring.  Indeed, both $I + \langle x_{\ell}
\rangle$ and $(I : x_{\ell})$ are strictly larger ideals than $I$, so
non-termination would give an infinite chain of strictly increasing
ideals.
\end{proof}

\begin{remark} \label{r:tree}
Running Algorithm~\ref{a:comma-colon} generates a rooted binary tree.
The nodes are monomial ideals and the root is the input ideal.  At
each node, Step~2 chooses a variable $x_{\ell}$.  The left-hand child
of a node $J$ is the ideal $J + \langle x_{\ell} \rangle$ and the
right-hand child is $(J : x_{\ell})$.  The corresponding branches are
labeled with the monomials $1$ and $x_{\ell}$ respectively.  The
leaves of this tree are prime ideals and each leaf corresponds to an
element in the Stanley decomposition.  Specifically, a leaf
corresponds to the pair $(\bm{x}^{\bm{u}}, \sigma)$ where
$\bm{x}^{\bm{u}}$ is the product of labels in the path from the root
to the leaf and $\sigma$ corresponds to the variables not in the prime
ideal.  We will call such a tree the associated binary tree for the
Stanley decomposition.  These trees also appear in \cite{clean}.
\end{remark}

\begin{example} \label{e:sdeg}
If $I = \langle x_{1}^{2}x_{2}, x_{1}x_{2}x_{3}, x_{2}^{2}x_{3},
x_{1}^{2}x_{4}, x_{1}x_{2}x_{4}, x_{2}^{2}x_{4} \rangle \subset S =
\Bbbk[x_{1}, x_{2}, x_{3}, x_{4}]$, then the Stanley decomposition
\[
\bigl\{ \bigl( 1, \{ 3, 4 \} \bigr), \bigl( x_{2}, \{ 3, 4 \} \bigr),
\bigl( x_{2}^{2} , \{ 2 \} \bigr), \bigl( x_{1} , \{ 3, 4 \} \bigr),
\bigl( x_{1}x_{2} , \{ 2 \} \bigr), \bigl( x_{1}^{2} , \{ 1, 3 \}
\bigr) \bigr\}
\] 
for $S/I$ produced by Algorithm~\ref{a:comma-colon} corresponds to the
following binary tree.
\[
\xymatrix@R=1.7em@C=-0.5mm{ 
& & & I \ar@{-}[dll]_{1} \ar@{-}[drr]^{x_{1}} \\
& \text{{\footnotesize $\langle x_{1}, x_{2}^{2}x_{3}, x_{2}^{2}x_{4}
\rangle$}}
\ar@{-}[dl]_{1} \ar@{-}[dr]^{x_{2}} & & & & \text{{\footnotesize
$\langle x_{1}x_{2}, x_{2}x_{3}, x_{1}x_{4}, x_{2}x_{4} \rangle$}}
\ar@{-}[dl]_{1} \ar@{-}[dr]^{x_{1}} \\
\text{{\footnotesize $\langle x_{1}, x_{2} \rangle$}} & &
\text{{\footnotesize $\langle x_{1}, x_{2}x_{3}, x_{2}x_{4} \rangle$}}
\ar@{-}[dl]_{1} \ar@{-}[d]^{x_{2}} & & \text{{\footnotesize $\langle
x_{1}, x_{2}x_{3}, x_{2}x_{4} \rangle$}} \ar@{-}[d]_{1}
\ar@{-}[dr]^{x_{2}} & & \text{{\footnotesize $\langle x_{2}, x_{4}
\rangle$}} \\
& \text{{\footnotesize $\langle x_{1}, x_{2} \rangle$}} &
\text{{\footnotesize $\langle x_{1}, x_{3}, x_{4} \rangle$}} & &
\text{{\footnotesize $\langle x_{1}, x_{2} \rangle$}} &
\text{{\footnotesize $\langle x_{1}, x_{3}, x_{4} \rangle$}} }
\]
\end{example}

These binary trees equip the Stanley decompositions produced by
Algorithm~\ref{a:comma-colon} with an additional structure.  To
describe this structure, we introduce the following concept.

\begin{definition}
A \emph{Stanley filtration} is a Stanley decomposition with an
ordering of the pairs $\bigl\{ (\bm{x}^{\bm{u}_{i}}, \sigma_{i}) : 1
\leq i \leq m \bigr\}$ such that for all $1 \leq j \leq m$ the set
$\bigl\{ (\bm{x}^{\bm{u}_{i}}, \sigma_{i}) : 1 \leq i \leq j \bigr\}$
is a Stanley decomposition for $S/\bigl(I + \langle
\bm{x}^{\bm{u}_{j+1}}, \dotsc, \bm{x}^{\bm{u}_{m}} \rangle \bigr)$.
Equivalently, the ordered set is a Stanley filtration provided the
modules $M_{j} = S/\bigl(I + \langle \bm{x}^{\bm{u}_{j+1}}, \dotsc,
\bm{x}^{\bm{u}_{m}} \rangle \bigr)$ form a filtration $\Bbbk = M_{0}
\subset M_{1} \subset \dotsc \subset M_{m} = S/I$ with $M_{j} /
M_{j-1} \cong S_{\sigma_{j}}$.
\end{definition}

\begin{example}
Not every Stanley decomposition has an ordering that makes it a
Stanley filtration.  For example, no ordering of the pairs in the
Stanley decomposition
\[
\bigl\{ \bigl( 1,\emptyset \bigr), \bigl(
x_{1}, \{ 1,2 \} \bigr), \bigl( x_{2}, \{ 2, 3 \} \bigr), \bigl(
x_{3}, \{ 1, 3 \} \bigr) \bigr\}
\] 
for $\Bbbk[x_{1}, x_{2}, x_{3}] / \langle x_{1}x_{2}x_{3} \rangle$ is
a Stanley filtration.
\end{example}

If $S/I$ has a Stanley filtration in which each $\sigma_{i}$
corresponds to a minimal prime of the ideal $I$, then Corollary~2.2.4
in \cite{clean} implies that $S/I$ is Cohen-Macaulay.

A standard way to traverse the leaves of a rooted tree is via
\emph{depth-first search} where all left-hand descendants of a node
are listed before any right-hand descendants.  This corresponds to
listing the leaves from left to right in the diagram of
Example~\ref{e:sdeg}.

\begin{corollary} \label{c:filtration}
Let $I \subseteq S$ be a monomial ideal and let $\mathfrak{S}$ be a
Stanley decomposition for $S/I$ obtained by applying
Algorithm~\ref{a:comma-colon}.  If the pairs have the order induced by
a depth-first search (starting with left-hand children) of the
associated binary tree, then $\mathfrak{S}$ is a Stanley filtration.
\end{corollary}

\begin{proof}
Let $\mathfrak{S} = \bigl\{ (\bm{x}^{\bm{u}_{i}}, \sigma_{i}) : 1 \leq
i \leq m \bigr\}$ and let $T$ be the binary tree associated to
$\mathfrak{S}$.  We write $L_{i}$ for the leaf corresponding to the
pair $(\bm{x}^{\bm{u}_{i}}, \sigma_{i})$.  We assume that $i < j$
implies that a depth-first search of $T$ arrives at $L_{i}$ before
reaching $L_{j}$.  It suffices to show that the set $\bigl\{
(\bm{x}^{\bm{u}_{i}}, \sigma_{i}) : 1 \leq i \leq m-1 \bigr\}$ can be
obtained by applying Algorithm~\ref{a:comma-colon} to $I + \langle
\bm{x}^{\bm{u}_{m}} \rangle$.  To accomplish this, we describe the
binary tree $T'$ generated by applying Algorithm~\ref{a:comma-colon}
to $I + \langle \bm{x}^{\bm{u}_{m}} \rangle$.  The tree $T'$ is
obtained from $T$ by deleting $L_{m}$ and contracting the branch
joining the parent of $L_{m}$ with its left-hand child.  The only
nodes in $T'$ that differ from $T$ are the first $| \bm{u}_{m} |$
nodes on the extreme right-hand branch.  These ideals are obtained
from those in $T$ by adding a proper divisor of $\bm{x}^{\bm{u}_{m}}$.
\end{proof}

\begin{example}
The converse of Corollary~\ref{c:filtration} is false, as there are
Stanley filtrations that do not arise from
Algorithm~\ref{a:comma-colon}.  For example, the third Stanley
decomposition in Example~\ref{e:sdegs} is a Stanley filtration with
respect to the given ordering that cannot be obtained from
Algorithm~\ref{a:comma-colon}.  Indeed, any decomposition obtained
from Algorithm~\ref{a:comma-colon} must have a term $( 1, \{ 1 \} )$
or $( 1, \{ 2 \})$ because $I + \langle x_{i} \rangle = \langle x_{i}
\rangle$ for $i = 1$, $2$.
\end{example}

%%---------------------------------------------------------------
\section{Bounds on Regularity} \label{s:boundandalgo}

This section contains the main results of this paper.  We first show
how a Stanley filtration for $S/I$ leads to a bound on its multigraded
regularity.

\begin{theorem} \label{t:sdbound}
Let $I$ be a monomial ideal in $S$.  If $\bigl\{ (
\bm{x}^{\bm{u}_{i}}, \sigma_{i}) : 1 \leq i \leq m \bigr\}$ is a
Stanley filtration for $S/I$, then $\bigcap_{i = 1}^{m} \bigl( A
\bm{u}_{i} + \reg(S_{\sigma_{i}}) \bigr) \subseteq \reg \left( S/I
\right)$.  In addition, if $I$ is $B$-saturated, then the intersection
can be taken over those pairs $(\bm{x}^{\bm{u}_{i}}, \sigma_{i})$ such
that $\widehat{\sigma}_{i} \in \Delta$.
\end{theorem}

\begin{proof}
Let $\scr{R}^{0}(M) := \{ \bm{k} \in \mathbb{Z}^{r} :
\text{$H_{B}^{0}(M)_{\bm{k} + \bm{c}} = 0$ for all $\bm{0} \neq \bm{c}
\in \scr{K}$} \}$ and for $j > 0$ set $\scr{R}^{j}(M) := \{ \bm{k} \in
\mathbb{Z}^{r} : \text{ $H_{B}^{j}(M)_{\bm{k} - \lambda_{1} \bm{c}_{1}
- \dotsb - \lambda_{e} \bm{c}_{e}} = 0$ for all $\lambda_{i} \in
\mathbb{N}$ with $\sum \lambda_{i} = j-1$} \}$.  With this notation,
we have $\reg(S/I) = \bigcap_{j \geq 0} \scr{R}^{j}(S/I)$.  We claim
that
\begin{equation} \label{sdbound:1}
\scr{R}^{j}(S/I) = \bigcap_{i = 1}^{m} \bigl( A \bm{u}_{i} +
\scr{R}^{j}(S_{\sigma_{i}}) \bigr) \, .
\end{equation}
This implies the first part of the theorem.  Additionally, if $I$ is
$B$-saturated then $\scr{R}^{0}(S/I) = \mathbb{Z}^{r}$ and $\reg(S/I)
= \bigcap_{j > 0} \scr{R}^{j}(S/I)$.  When $\widehat{\sigma} \not\in
\Delta$, Lemma~\ref{l:Bsaturated} implies that $S_{\sigma} =
S/P_{\sigma}$ is a $B$-torsion module, so $H_{B}^{j}(S_{\sigma}) = 0$
for $j > 0$.  It follows that $\scr{R}^{j}(S_{\sigma}) =
\mathbb{Z}^{r}$ for $j > 0$ and hence $\bigcap_{i = 1}^{m} \bigl( A
\bm{u}_{i} + \scr{R}^{j}(S_{\sigma_{i}}) \bigr) =
\bigcap_{\widehat{\sigma}_{i} \in \Delta} \bigl( A \bm{u}_{i} +
\scr{R}^{j}(S_{\sigma_{i}}) \bigr)$.  Therefore, the claim also
establishes the second part of the theorem.

We prove \eqref{sdbound:1} by induction on $m$.  When $m = 1$, the
unique pair has the form $(1, \sigma)$ which implies that $I =
P_{\sigma} = \langle x_{i} : i \not\in \sigma \rangle$ and
$\scr{R}^{j}(S/I) = \scr{R}^{i}(S_{\sigma})$.  Suppose that the claim
holds for all Stanley filtrations with fewer than $m$ pairs.  The
short exact sequence $0 \longrightarrow S(-A \bm{u}_{m})/(I :
\bm{x}^{\bm{u}_{m}}) \xrightarrow{\; \bm{x}^{\bm{u}_{m}} \;} {S}/{I}
\longrightarrow {S}/(I+ \langle \bm{x}^{\bm{u}_{m}} \rangle)
\longrightarrow 0$ yields the exact sequence
\begin{equation} \label{longexact}
H_{B}^{j} \bigl( S/(I : \bm{x}^{\bm{u}_{m}}) \bigr)_{\bm{p} - A
\bm{u}_{m}} \longrightarrow H_{B}^{j} \bigl( S/I \bigr)_{\bm{p}}
\longrightarrow H_{B}^{j} \bigl( S/(I + \langle \bm{x}^{\bm{u}_{m}}
\rangle) \bigr)_{\bm{p}} \, .
\end{equation}
From this, we deduce that $\scr{R}^{j} \bigl( S/(I + \langle
\bm{x}^{\bm{u}_{m}}) \bigr) \cap \bigl( A \bm{u}_{m} + \scr{R}^{j}
\bigl( S / (I : \bm{x}^{\bm{u}_{m}}) \bigr) \bigr) \subseteq
\scr{R}^{j}(S/I)$.  Since $\bigl\{ (\bm{x}^{\bm{u}_{i}}, \sigma_{i}) :
1 \leq i \leq m-1 \bigr\}$ is a Stanley filtration for $S/(I + \langle
\bm{x}^{\bm{u}_{m}} \rangle)$, the induction hypothesis implies that
$\scr{R}^{j} \bigl( S / (I + \langle \bm{x}^{\bm{u}_{m}} \rangle)
\bigr) = \bigcap_{i=1}^{m-1} \bigl( A \bm{u}_{i} +
\scr{R}^{j}(S_{\sigma_{i}}) \bigr)$.  The ordering also implies that
no monomial in $S$ divisible by $\bm{x}^{\bm{u}_{m}}$ belongs to the
set $\bigcup_{i=1}^{m-1} (\bm{x}^{\bm{u}_{i}}, \sigma_{i})$.  It
follows that a monomial $\bm{x}^{\bm{u}_{m} + \bm{v}} \in S$ is not
contained in $I$ if and only if $\supp(\bm{v}) \subseteq \sigma_{m}$.
Therefore, we have $(I : \bm{x}^{\bm{u}_{m}}) = P_{\sigma_{m}}$ and
$\scr{R}^{j} \bigl(S/(I : \bm{x}^{\bm{u}_{m}}) \bigr) =
\scr{R}^{j}(S_{\sigma_{m}})$ which completes the induction.
\end{proof} 

\begin{remark} \label{r:totaldegree}
If $S$ has the standard grading (equivalently $X = \mathbb{P}^{d}$),
then Theorem~\ref{t:sdbound} says that the Castelnuovo-Mumford
regularity of a monomial ideal is bounded by the maximum of
$|\bm{u}_{i}| := \sum_{i = 1}^{n} u_{i}$ for a Stanley filtration
$\bigl\{ ( \bm{x}^{\bm{u}_{i}}, \sigma_{i}) \bigr\}$.
\end{remark}

We next examine the relationship between Stanley filtrations and
Hilbert polynomials.  Given a Stanley filtration $\bigl\{
(\bm{x}^{\bm{u}_{i}}, \sigma_{i}) : 1 \leq i \leq m \bigr\}$ for
$S/I$, we have
\[
H( S/I, \bm{t}) = \sum_{i=1}^{m} H( S_{\sigma_{i}}, \bm{t} - A
\bm{u}_{i}) \, .
\]  
Since $\scr{K} \subseteq \mathbb{N} A_{\sigma}$
if and only if $\widehat{\sigma} \in \Delta$, the Hilbert polynomial
of $S/I$ has an expression with potentially fewer summands:
$P_{S/I}(\bm{t}) = \sum_{\widehat{\sigma} \in \Delta}
P_{S_{\sigma}}(\bm{t} - A \bm{u})$.  To place further restrictions on
the summands, we need an ordering on the $\widehat{\sigma} \in
\Delta$.

We endow the polynomial ring $\mathbb{Q}[t_{1}, \dotsc, t_{r}]$ with
the $\mathbb{Z}$-grading defined by $\deg(t_{i}) = 1$ for $1 \leq i
\leq r$.  Let $<$ be a monomial order on $\mathbb{Q}[t_{1}, \dotsc,
t_{r}]$ which refines the order by total degree.  This graded monomial
order induces a partial order, also denoted $<$, on the simplices of
$\Delta$.  Specifically, $\widehat{\sigma} < \widehat{\tau}$ if and
only if $\initial_{<}\bigl( P_{S_{\tau}}(\bm{t}) \bigr) <
\initial_{<}\bigl( P_{S_{\sigma}}(\bm{t}) \bigr)$.  Since the total
degree of $P_{S_{\sigma}}(\bm{t})$ equals $|\sigma| - d$, the induced
order on $\Delta$ refines inclusion: $\widehat{\sigma} \subseteq
\widehat{\tau}$ implies $\widehat{\sigma} \leq \widehat{\tau}$.

\begin{definition}
A total order $\prec$ on $\Delta$ is called \emph{graded} if it
refines the partial order induced by a graded monomial order $<$ on
$\mathbb{Q}[t_{1}, \dotsc, t_{r}]$.
\end{definition}

\begin{proposition} \label{p:niceSD}
If $\prec$ is a graded total order on $\Delta$ and $I \subset S$ is a
monomial ideal, then $S/I$ has a Stanley filtration $\bigl\{
(\bm{x}^{\bm{u}_{i}}, \sigma_{i}) : 1 \leq i \leq m \bigr\}$
satisfying the following condition:
\begin{itemize}
\item if there is an index $i$ with $\widehat{\sigma}_{i} \in \Delta$
and $\bm{x}^{\bm{u}_{i}} \neq 1$, then there exists an
index $j < i$ such that $\widehat{\sigma}_{j} \in \Delta$,
$\widehat{\sigma}_{j} \preceq \widehat{\sigma}_{i}$ and
$\bm{x}^{\bm{u}_{i}} = \bm{x}^{\bm{u}_{j}}
x_{\ell}$ for some $\ell \not\in \sigma_{j}$.
\end{itemize}
\end{proposition}

\begin{proof}
We refine Step~$2$ of Algorithm~\ref{a:comma-colon} to produce a
Stanley filtration that satisfies the given condition.  Specifically,
Step~$2$ becomes:
\begin{enumerate}
\item[$2'$.] (Choose variable) If $I$ is not contained in $P_{\tau}$
for some $\widehat{\tau} \in \Delta$, then choose a variable
$x_{\ell} \in S$ that is a proper divisor of a minimal generator of
$I$.  Otherwise, let $\widehat{\sigma} \in \Delta$ be the smallest
simplex with respect to $\prec$ for which $I \subsetneq P_{\sigma} =
\langle x_{i} : i \in \widehat{\sigma} \rangle$ and choose a variable
$x_{\ell} \in P_{\sigma}$ that is a proper divisor of a minimal
generator of $I$.
\end{enumerate}
To prove that the resulting Stanley filtration has the desired form,
we analyze the associated binary tree.  Let $(\bm{x}^{\bm{u}_{i}},
\sigma_{i})$ be a pair in the Stanley filtration with
$\widehat{\sigma}_{i} \in \Delta$ and let $L_{i}$ be the
corresponding leaf.  The leaf $L_{i}$ is either a left-hand or
right-hand child of its parent.

Suppose $L_{i}$ is a right-hand child.  We write $J$ for the parent of
$L_{i}$ and $x_{\ell}$ for the variable labeling the branch connecting
$J$ and $L_i$, so $(J : x_{\ell}) = L_{i}$.  Let $L_{j}$ be the
descendant of $J$ obtained by repeatedly taking the left-hand child of
$J$.  The leaf $L_{j}$ corresponds to a pair $(\bm{x}^{\bm{u}_{j}},
\sigma_{j})$.  Since the left-hand branches are always labeled with
$1$, we see that $\bm{x}^{\bm{u}_{j}} x_{\ell} = \bm{x}^{\bm{u}_{i}}$.
Moreover, the depth-first search ordering (see
Corollary~\ref{c:filtration}) chooses left-hand children first, so we
have $j < i$.  Because all the left-hand descendants of $J$ contain
$x_{\ell}$, we must also have $\ell \not\in \sigma_{j}$.

It remains to show that $\widehat{\sigma}_{j} \in \Delta$ and
$\widehat{\sigma}_{j} \preceq \widehat{\sigma}_{i}$.  Because $L_{i}
= P_{\sigma_{i}}$, we have $J \subset P_{\sigma_{i}}$.  Hence the set
of all $\widehat{\sigma} \in \Delta$ with $J \subset P_{\sigma}$ is
nonempty, so we may take $\widehat{\tau} \in \Delta$ to be one which
is minimal with respect to $\prec$.  Step~$2'$ guarantees that every
left-hand child of $J$ is also contained in $P_{\tau}$.  This
containment must be proper until the leaf $L_{j}$ is reached.  This
means that $P_{\sigma_{j}} \subseteq P_{\tau}$ which implies
$\widehat{\sigma}_{j} \subseteq \widehat{\tau}$ and
$\widehat{\sigma}_{j} \in \Delta$.  Since $J \subset P_{\sigma_{i}}$,
the minimality of $\widehat{\tau}$ implies that $\widehat{\tau}
\preceq \widehat{\sigma}_{i}$.  Hence, we have $\widehat{\sigma}_{j}
\in \Delta$ and $\widehat{\sigma}_{j} \preceq \widehat{\sigma}_{i}$
as required.

On the other hand, suppose that $L_{i}$ is a left-hand child.  Let
$J'$ be the closest ancestor of $L_{i}$ that is a right-hand child.
Such an ancestor exists if and only if $\bm{x}^{\bm{u}_{i}} \neq 1$.
Since $J' \subset P_{\sigma_{i}}$, the argument given when $L_{i}$ is
a right-hand child applies to the parent of $J'$ and this completes
the proof.
\end{proof}

Using Proposition~\ref{p:niceSD}, we can give an algorithm for finding
all $B$-saturated monomial ideals with a given Hilbert polynomial
$P(\bm{t})$.  Roughly speaking, the algorithm works by ``peeling off''
smaller Hilbert polynomials $P_{S_{\sigma}}(\bm{t})$ from $P(\bm{t})$.
To accomplish this, we need the following result about the leading
coefficients of the Hilbert polynomial.  This lemma generalizes
techniques used in the proof of Theorem~3.2 of \cite{HoangTrung}.

\begin{lemma} \label{l:leading}
Let $\bm{e}_{1} := [ 1 \; 0 \; \dotsb \; 0 ]^{\text{\sf{T}}} \in
\mathbb{Z}^{r}$ be the first standard basis vector and let $P(\bm{t})$
be the multigraded Hilbert polynomial of $M$.  If $\bm{e}_{1} \in
\intt \scr{K}$, then the leading coefficient of $P(\bm{t})$ with
respect to the graded lexicographic order with $t_{1} > t_{2} > \dotsb
> t_{r}$ is positive.
\end{lemma}

\begin{proof}
Using Proposition~1.11 in \cite{sturmfels}, we may choose a weight
vector $\bm{w} \in \mathbb{N}^{r}$ such that $\initial_{\bm{w}}(P) =
\initial_{\text{glex}}(P)$ and $w_{1} > w_{i}$ for $1 < i \leq r$.
Let $\varphi_{\bm{w}} \colon \mathbb{N} \to \mathbb{N}^{r}$ be the map
defined by $\varphi_{\bm{w}}(z) = (z^{w_{1}}, \dotsc, z^{w_{r}})$.
Since $w_{1}$ is the largest component of $\bm{w}$, we have
$\lim\limits_{z \to \infty} \tfrac{\varphi_{\bm{w}}(z)}{\|
\varphi_{\bm{w}}(z) \|} = \bm{e}_{1}$.  By hypothesis, we also have
$\bm{e}_{1} \in \intt \scr{K}$ which implies that $\varphi_{\bm{w}}(z)
\in \scr{K}$ for $z \gg 0$.  For a fixed $z$, consider $Q_{z}(y) =
P\bigl(y \, \varphi_{\bm{w}}(z) \bigr) \in \mathbb{Q}[y]$.  If
$P(\bm{t}) = \sum_{\bm{u}} b_{\bm{u}} \, \bm{t}^{\bm{u}}$ has total
degree $\ell$ then the highest degree term in $Q_{z}(y)$ is $\left(
\sum_{| \bm{u}| = \ell} b_{\bm{u}} z^{\bm{w} \cdot \bm{u}} \right)
y^{\ell}$.  When $\varphi_{\bm{w}}(z) \in \scr{K}$ and $y \gg 0$,
$Q_{z}(y)$ agrees with the Hilbert function $H \bigr( M, y \,
\varphi_{\bm{w}}(z) \bigl)$ which implies that $Q_{z}(y) > 0$.  Thus,
the leading coefficient of the polynomial $Q_{z}(y)$ is positive.
Because this is true for all sufficiently large $z$, the leading
coefficient of $\sum_{|\bm{u}| = \ell} b_{\bm{u}} z^{\bm{w} \cdot
\bm{u}}$ considered as a polynomial in ${\mathbb Q}[z]$ is also
positive.  Finally, our choice of $\bm{w}$ implies that the leading
coefficient of $\sum_{|\bm{u}| = \ell} b_{\bm{u}} z^{\bm{w} \cdot
\bm{u}}$ equals the leading coefficient of $P(\bm{t})$ with respect to
the graded lexicographic order.
\end{proof}

\begin{remark}
Proposition~\ref{l:leading} is more applicable than is obvious at
first glance.  Clearly $\bm{e}_1$ can be replaced by any other
standard basis vector $\bm{e}_i$, with the corresponding change of
lexicographic order.  More generally, there is always a unimodular
coordinate change on $\mathbb{Z}^{r}$ that takes the configuration $\{
\bm{a}_{1}, \dotsc, \bm{a}_{n} \}$ to a new configuration $\{
\bm{a}_{1}', \dotsc, \bm{a}_{n}' \}$ with $\bm{e}_{1} \in \intt
\scr{K}$.  Indeed, any vector $\bm{v} \in \mathbb{Z}^{r}$ with
$\gcd(v_{i}) = 1$ can be the first column of a matrix in
$\text{SL}_{r}(\mathbb{Z})$.  In fact, there is an unimodular
transformation of $\mathbb{Z}^{r}$ such that the entire positive
orthant lies inside $\scr{K}$.  In this case, the leading term of the
Hilbert polynomial $P(\bm{t})$ with respect to any graded monomial
order (not just graded lexicographic ones) on $\Bbbk[t_{1}, \dotsc,
t_{r}]$ is positive.  This conclusion also holds provided the sequence
$\varphi_{\bm{w}}(z)$ approaches $\| \varphi_{\bm{w}}(z) \|
\bm{e}_{i}$ from within $\scr{K}$.  In particular, it applies when
$\scr{K}$ equals the positive orthant as in
Examples~\ref{e:examplecones} and \ref{e:hilbertpoly}.
\end{remark}

We now use Proposition~\ref{p:niceSD} and Lemma~\ref{l:leading} to
give an algorithm for listing all $B$\nobreakdash-saturated monomial
ideals with a given multigraded Hilbert polynomial.

\begin{algorithm} \label{a:macaulay}
Let $\prec$ be a graded total order on $\Delta$ and let $<$ be the
corresponding graded monomial order on $\mathbb{Q}[t_{1}, \dotsc,
t_{r}]$.  Given a polynomial $P(\bm{t}) \in \mathbb{Q}[t_{1}, \dotsc,
t_{r}]$, this algorithm returns all $B$-saturated monomial ideals with
the multigraded Hilbert polynomial $P(\bm{t})$.
\begin{enumerate}
\item (Coordinate change) If necessary, make a unimodular coordinate
change $\phi$ on $\mathbb{Z}^{r}$ such that the positive orthant lies
inside $\scr{K}$ and replace the polynomial $P(\bm{t})$ with $P
\bigl( \phi^{-1}(\bm{t}) \bigr)$.
\item (Initialize) Set $\textsf{Reps} = \emptyset$,
$\textsf{PartialReps} = \bigl\{ \bigl( \emptyset, P(\bm{t}) \bigr)
\bigr\}$ and $\textsf{Ideals} = \emptyset$.
\item (Enlarge representation) Select and remove an element $\bigr(
\mathfrak{S}, Q(\bm{t}) \bigr) \in \textsf{PartialReps}$.  For each
$\widehat{\tau} \in \Delta$ satisfying
\begin{itemize}
\item[(a)] if $\mathfrak{S} \neq \emptyset$, then there exists a pair
$(\bm{x}^{\bm{u}}, \sigma) \in \mathfrak{S}$ with $\widehat{\sigma}
\preceq \widehat{\tau}$;
\item[(b)] $\initial_{<} \bigl( Q(\bm{t}) \bigr) = \initial_{<} \bigl(
P_{S_{\tau}}(\bm{t}) \bigr)$;
\item[(c)] the leading coefficient with respect to $<$ of $Q(\bm{t}) -
P_{S_{\tau}}(\bm{t})$ is positive;
\end{itemize}
and for each monomial $\bm{x}^{\bm{v}} \in S$ satisfying
\begin{itemize}
\item[(d)] if $\mathfrak{S} = \emptyset$ then $\bm{x}^{\bm{v}} =
1$;
\item[(e)] if $\mathfrak{S} \neq \emptyset$ then for
$(\bm{x}^{\bm{u}}, \sigma)$ from (a) we have $\bm{x}^{\bm{v}} =
\bm{x}^{\bm{u}} x_{\ell}$ for some $\ell \not \in \sigma$;
\end{itemize}
do as follows.  If $Q(\bm{t}) - P_{S_{\tau}}(\bm{t}) = 0$ then append
the set $\mathfrak{S} \cup \{ (\bm{x}^{\bm{v}}, \tau) \}$ to
$\textsf{Reps}$.  Otherwise, append the pair $\bigl( \mathfrak{S} \cup
\{ (\bm{x}^{\bm{v}}, \tau) \}, Q(\bm{t}) - P_{S_{\tau}}(\bm{t})
\bigr)$ to $\textsf{PartialReps}$.
\item (Finished?)  If $\textsf{PartialReps} \neq \emptyset$ then go to
step~3.
\item (Check Hilbert polynomial) For each $\mathfrak{S} \in
\textsf{Reps}$ compute the multigraded Hilbert polynomial of the ideal
$I = \bigcap\nolimits_{(\bm{x}^{\bm{v}},\tau) \in \mathfrak{S}}
\langle x_i^{v_i+1} : i \not \in \tau \rangle$.  If the multigraded
Hilbert polynomial of $I$ is $P(\bm{t})$ then append $I$ to
$\textsf{Ideals}$.  Output the list $\textsf{Ideals}$.
\end{enumerate}
\end{algorithm}

\begin{proof}[Proof of Correctness]
By construction, the output is a list of monomial ideals with
multigraded Hilbert polynomial $P(\bm{t})$ that are $B$-saturated by
Lemma~\ref{l:Bsaturated}.  Conversely, given any
$B$\nobreakdash-saturated monomial ideal $I$,
Proposition~\ref{p:niceSD} provides a Stanley filtration $\{(
\bm{x}^{\bm{u}_{i}}, \sigma_{i}) : 1 \leq i \leq m \}$ for $S/I$ such
that for all $i > 1$ there is a $j < i$ with $\widehat{\sigma}_{j} \in
\Delta$, $\widehat{\sigma}_{j} \leq \widehat{\sigma}_{i}$ and
$\bm{x}^{\bm{u}_{i}} = \bm{x}^{\bm{u}_{j}}x_{\ell}$ for some $\ell
\not \in \sigma_j$.  Thus, the conditions~(a), (d) and (e) in Step~3
do not eliminate any $B$-saturated monomial ideals with Hilbert
polynomial $P(\bm{t})$.

For $1 \leq j \leq m$, the polynomial $P(\bm{t}) - \sum_{i=1}^{j}
P_{S_{\sigma_{i}}}(\bm{t})$ is the multigraded Hilbert polynomial of
the $\mathbb{Z}^{r}\!$-graded $S$-module $\bigoplus_{i=j+1}^{m}
S_{\sigma_{i}}(-A \bm{u}_{i})$ and Lemma~\ref{l:leading} (combined
with Step~1) ensures that its leading coefficient is positive.  Since
$\prec$ is a graded total order on $\Delta$, we have $\initial \bigl(
P_{S_{\sigma_{i}}}(\bm{t}) \bigr) \geq \initial \bigl(
P_{S_{\sigma_{j}}}(\bm{t}) \bigr)$ for $i < j$, so the leading term of
the subtracted polynomial will be $\initial(P_{S_{\sigma_{j+1}}})$.
This means that conditions~(b) and (c) in Step~3 do not exclude any of
the relevant ideals.  We conclude that every $B$-saturated monomial
ideal with multigraded Hilbert polynomial $P(\bm{t})$ has a Stanley
filtration of the form created by this procedure which implies every
such ideal is part of the output.

It remains to show that this procedure terminates.  To accomplish
this, observe that Step~3 replaces the pair $(\mathfrak{S}, Q)$ with
pairs in which either the leading coefficient of the second entry, or
its leading term, is strictly less than that of $Q(\bm{t})$.  Since
there are only a finite number of choices for $P_{S_{\tau}}(\bm{t})$,
there is a lower bound on how much the leading coefficient can
decrease which guarantees that the process cannot continue
indefinitely.
\end{proof}

This corollary also follows, albeit non-constructively, from
\cite{antichains}.
 
\begin{corollary} \label{c:finite}
For any polynomial $P(\bm{t})$, there are only finitely many
$B$-saturated monomial ideals with multigraded Hilbert polynomial
$P(\bm{t})$. \qed
\end{corollary}

We illustrate Algorithm~\ref{a:macaulay} by constructing all
$B$-saturated monomial ideals in the standard graded polynomial ring
$S = \Bbbk[x_{1}, x_{2}, x_{3}]$ having Hilbert polynomial $3t+1$.

\begin{example}
Since the lead term of the Hilbert polynomial is $3t$, there must be
three pairs of the form $(\bm{x}^{\bm{u}}, \tau)$ with $| \tau | = 2$.
Fix the ordering: $\{ 1 \} \prec \{ 2 \} \prec \{ 3 \}$.  Since the
pairs correspond to disjoint sets of monomials, the first three pairs
are $(1, \{i_{1}, i_{2} \})$, $(x_{i_{3}}, \{ j_{1}, j_{2} \})$ and
$(x_{i_{3}}x_{j_{3}}, \{ k_{1}, k_{2} \})$ where $\{i_{1}, i_{2},
i_{3} \} = \{j_{1}, j_{2}, j_{3} \} = \{1,2,3\}$.  These pairs
contribute $\binom{t+1}{1} + \binom{t}{1} + \binom{t-1}{1} = 3t$ to
the Hilbert polynomial.  Hence, the Stanley filtrations also contain
the pair $(x_{i_{3}}x_{j_{3}}x_{k_{3}}, \{ \ell_{1} \})$ where $\{
k_{1}, k_{2}, k_{3} \} = \{ 1, 2, 3\}$.  Constructing all these sets
which satisfy the order condition gives:
\begin{flalign*}
& \text{\scriptsize $\bigl\{ (1,\{ 2, 3 \}), (x_{1},\{ 2, 3 \}),
(x_{1}^{2},\{ 2, 3 \}), (x_{1}^{3},\{ 3 \}) \bigr\},$} & &
\text{\scriptsize $\bigl\{ (1,\{ 2, 3 \}), (x_{1},\{ 2, 3 \}),
(x_{1}^{2},\{ 2, 3 \}), (x_{1}^{3},\{ 2 \}) \bigr\},$} \\
& \text{\scriptsize $\bigl\{ (1,\{ 2, 3 \}), (x_{1},\{ 2, 3 \}),
(x_{1}^{2},\{ 2, 3 \}), (x_{1}^{3},\{ 1 \}) \bigr\},$} & &
\text{\scriptsize $\bigl\{ (1,\{ 2, 3 \}), (x_{1},\{ 2, 3 \}),
(x_{1}^{2},\{ 2, 3 \}), (x_{1}^{2}x_{2},\{ 3 \}) \bigr\},$} \\
& \text{\scriptsize $\bigl\{ (1,\{ 2, 3 \}), (x_{1},\{ 2, 3 \}),
(x_{1}^{2},\{ 1, 3 \}), (x_{1}^{2}x_{2}, \{ 2 \}) \bigr\},$} & &
\text{\scriptsize $\bigl\{ (1,\{ 2, 3 \}), (x_{1},\{ 2, 3 \}),
(x_{1}^{2},\{ 2, 3 \}), (x_{1}^{2}x_{2},\{ 1 \}) \bigr\},$} \\
& \text{\scriptsize $\bigl\{ (1,\{ 2, 3 \}), (x_{1},\{ 2, 3 \}),
(x_{1}^{2},\{ 1, 2 \}), (x_{1}^{2}x_{3}, \{ 2 \}) \bigr\},$} & &
\text{\scriptsize $\bigl\{ (1,\{ 2, 3 \}), (x_{1},\{ 2, 3 \}),
(x_{1}^{2},\{ 1, 2 \}), (x_{1}^{2}x_{3},\{ 1 \}) \bigr\},$} \\
& \text{\scriptsize $\bigl\{ (1,\{ 2, 3 \}), (x_{1},\{ 2, 3 \}),
(x_{1}^{2},\{ 1, 2 \}), (x_{1}^{2}x_{3}, \{ 3 \}) \bigr\},$} & &
\text{\scriptsize $\bigl\{ (1,\{ 2, 3 \}), (x_{1},\{ 1, 3 \}),
(x_{1}x_{2},\{ 1, 3 \}), (x_{1}x_{2}^{2}, \{ 3 \}) \bigr\},$} \\
& \text{\scriptsize $\bigl\{ (1,\{ 2, 3 \}), (x_{1},\{ 1, 3 \}),
(x_{1}x_{2},\{ 1, 3 \}), (x_{1}x_{2}^{2},\{ 2 \}) \bigr\},$} & &
\text{\scriptsize $\bigl\{ (1,\{ 2, 3 \}), (x_{1},\{ 1, 3 \}),
(x_{1}x_{2},\{ 1, 3 \}), (x_{1}x_{2}^{2}, \{ 1 \}) \bigr\},$} \\
& \text{\scriptsize $\bigl\{ (1,\{ 2, 3 \}), (x_{1},\{ 1, 3 \}),
(x_{1}x_{2},\{ 1, 2 \}), (x_{1}x_{2}x_{3}, \{ 3 \}) \bigr\},$} & &
\text{\scriptsize $\bigl\{ (1,\{ 2, 3 \}), (x_{1},\{ 1, 3 \}),
(x_{1}x_{2},\{ 1, 2 \}), (x_{1}x_{2}x_{3}, \{ 2 \}) \bigr\},$} \\
& \text{\scriptsize $\bigl\{ (1,\{ 2, 3 \}), (x_{1},\{ 1, 3 \}),
(x_{1}x_{2},\{ 1, 2 \}), (x_{1}x_{2}x_{3}, \{ 1 \}) \bigr\},$} & &
\text{\scriptsize $\bigl\{ (1,\{ 2, 3 \}), (x_{1},\{ 1, 2 \}),
(x_{1}x_{3},\{ 1, 2 \}), (x_{1}x_{3}^{2}, \{ 2 \}) \bigr\},$} \\
& \text{\scriptsize $\bigl\{ (1,\{ 2, 3 \}), (x_{1},\{ 1, 2 \}),
(x_{1}x_{3},\{ 1, 2 \}), (x_{1}x_{3}^{2}, \{ 1 \}) \bigr\},$} & &
\text{\scriptsize $\bigl\{ (1,\{ 1, 3 \}), (x_{2},\{ 1, 3 \}),
(x_{2}^{2}, \{ 1, 3 \}), (x_{2}^{3}, \{ 3 \}) \bigr\},$} \\
& \text{\scriptsize $\bigl\{ (1,\{ 1, 3 \}), (x_{2},\{ 1, 3 \}),
(x_{2}^{2}, \{ 1, 3 \}), (x_{2}^{3}, \{ 2 \}) \bigr\},$} & &
\text{\scriptsize $\bigl\{ (1,\{ 1, 3 \}), (x_{2},\{ 1, 3 \}),
(x_{2}^{2}, \{ 1, 3 \}), (x_{2}^{3}, \{ 1 \}) \bigr\},$} \\
& \text{\scriptsize $\bigl\{ (1,\{ 1, 3 \}), (x_{2},\{ 1, 3 \}),
(x_{2}^{2}, \{ 1, 2 \}), (x_{2}^{2}x_{3}, \{ 2 \}) \bigr\},$} & &
\text{\scriptsize $\bigl\{ (1,\{ 1, 3 \}), (x_{2},\{ 1, 3 \}),
(x_{2}^{2}, \{ 1, 2 \}), (x_{2}^{2}x_{3}, \{ 1 \}) \bigr\},$} \\
& \text{\scriptsize $\bigl\{ (1,\{ 1, 3 \}), (x_{2},\{ 1, 2 \}),
(x_{2}x_{3}, \{ 1, 2 \}), (x_{2}x_{3}^{2}, \{ 2 \}) \bigr\},$} & &
\text{\scriptsize $\bigl\{ (1,\{ 1, 3 \}), (x_{2},\{ 1, 2 \}),
(x_{2}x_{3}, \{ 1, 2 \}), (x_{2}x_{3}^{2}, \{ 1 \}) \bigr\},$} \\
& \text{\scriptsize $\bigl\{ (1,\{ 1, 2 \}), (x_{3},\{ 1, 2 \}),
(x_{3}^{2}, \{ 2, 3 \}), (x_{1}x_{3}^{2}, \{ 3 \}) \bigr\},$} & &
\text{\scriptsize $\bigl\{ (1,\{ 1, 2 \}), (x_{3},\{ 1, 2 \}),
(x_{3}^{2}, \{ 1, 3 \}), (x_{2}x_{3}^{2}, \{ 3 \}) \bigr\},$} \\
& \text{\scriptsize $\bigl\{ (1,\{ 1, 2 \}), (x_{3},\{ 1, 2 \}),
(x_{3}^{2}, \{ 1, 2 \}), (x_{3}^{3}, \{ 3 \}) \bigr\},$} & &
\text{\scriptsize $\bigl\{ (1, \{ 1, 2 \}), (x_{3},\{ 1, 3 \}),
(x_{2}x_{3}, \{ 1, 3 \}), (x_{2}^{2}x_{3}, \{ 3 \}) \bigr\},$} \\
& \text{\scriptsize $\bigl\{ (1,\{ 1, 2 \}), (x_{3},\{ 1, 2 \}),
(x_{3}^{2}, \{ 1, 2 \}), (x_{3}^{3}, \{ 2 \}) \bigr\},$} & &
\text{\scriptsize $\bigl\{ (1, \{ 1, 2 \}), (x_{3},\{ 1, 2 \}),
(x_{2}x_{3}, \{ 1, 2 \}), (x_{3}^{3}, \{ 1 \}) \bigr\}.$} 
\end{flalign*}
In particular, there are $30$ $B$-saturated monomial ideals in $S$
with Hilbert function $3t+1$.

We can verify this calculation as follows.  Since $3t+1 =
\binom{t+1}{1} + \binom{t}{1}+ \binom{t-1}{1} + \binom{t-2}{0}$,
Gotzmann's regularity theorem implies that every saturated ideal with
the required Hilbert polynomial has regularity $4$ which means the
generators have degree at most $4$.  Because $\dim_{\Bbbk} S_{\bm{4}}
= 15$ and $3(4) + 1 = 13$, the list consists of all ideals generated
by two monomials of degree $4$.  Eliminating those that do not have
the correct Hilbert polynomial produces the same $30$ monomial ideals.
\end{example}

To state our next theorem, we make the following definition.

\begin{definition} \label{d:gotzmannnum}
Let $m$ be the largest number of pairs in a decomposition
$\mathfrak{S}$ constructed in Algorithm~\ref{a:macaulay}.  We call
this the \emph{Gotzmann number} of $P(\bm{t})$.
\end{definition}

To calculate an upper bound for the Gotzmann number of $P(\bm{t})$, we
can use a simplified version of Algorithm~\ref{a:macaulay}.
Specifically, the Gotzmann number is bounded by the maximum $k$ among
all the expressions $P(\bm{t}) = \sum_{i = 1}^{k}
P_{i}(\bm{t}-\bm{q}_{i})$ that satisfy the following conditions:
\begin{enumerate}
\item $P_{i}(\bm{t}) = P_{S_{\sigma_{i}}}(\bm{t})$ for some
$\widehat{\sigma}_{i} \in \Delta$;
\item $\bm{q}_{1}=0$;
\item for all $i > 1$, there is a $j < i$ with $\widehat{\sigma_{j}}
\preceq \widehat{\sigma_{i}}$ and $\bm{q}_{i} = \bm{q}_{j} +
\bm{a}_{\ell}$ for some $\ell \not\in \sigma_{j}$.
\end{enumerate}
When $S$ has the standard grading (or equivalently when $X =
\mathbb{P}^{d}$), the results of \S\ref{s:standardcase} show that this
upper bound is the exact Gotzmann number.  The analogous question for
general smooth projective toric varieties is not known.

We now establish our multigraded analogue of Gotzmann's regularity
theorem.

\begin{theorem} \label{t:hilbpolybound}
Let $I$ be any $B$-saturated ideal in $S$ and let $\bm{c} \in
\bigcap_{i=1}^{n} \bigl( \bm{a}_{i} + \scr{K} \bigr)$.  If $m$ is the
Gotzmann number of the Hilbert polynomial $P_{S/I}(\bm{t})$ then
\[
\bigcap\limits_{\widehat{\sigma} \in \Delta} \bigl( (m-1) \bm{c} +
\reg(S_{\sigma}) \bigr) \subseteq \reg(S/I) \, .
\]
\end{theorem}

\begin{proof}
Applying Proposition~\ref{p:usc} and Lemma~\ref{l:saturationpoly}, we
may assume without loss of generality that $I$ is a $B$-saturated
monomial ideal.  Algorithm~\ref{a:macaulay} yields a partial Stanley
filtration $\bigl\{ (\bm{x}^{\bm{u}_{i}}, \sigma_{i}) \bigr\}$ with at
most $m$ pairs.  Moreover, we have $|\bm{u}_{i}| < i$.  Since the
hypothesis on $\bm{c}$ guarantees that $(m-1) \bm{c} +
\reg(S_{\sigma_{i}}) \subseteq A \bm{u}_{i} + \reg(S_{\sigma_{i}})$
and Theorem~\ref{t:sdbound} implies that
$\bigcap_{\widehat{\sigma}_{i} \in \Delta} \bigl( A \bm{u}_{i} +
\reg(S_{\sigma_{i}}) \bigr) \subseteq \reg \left( S/I \right)$, the
theorem follows.
\end{proof}

We end this section with two examples.

\begin{example}
Let $I$ be an $B$-saturated ideal corresponding to the set of $\ell$
points on a smooth projective toric variety $X$.  Hence,
$P_{S/I}(\bm{t}) = \ell$ and the Gotzmann number of $P_{S/I}(\bm{t})$
is also $\ell$.  If $\bm{c} \in \bigcap_{i = 1}^{n} (\bm{a}_{i} +
\scr{K})$, then $(\ell-1)\bm{c}$-regular.  This bound is independent
of the configuration of the points.  In contrast, Proposition~6.7 in
\cite{MaclaganSmith1} shows that $\reg(S/I)$ does depend on the
arrangement the points on $X$.
\end{example}

\begin{example}
If $X = \mathbb{P}^{2} \times \mathbb{P}^{1}$ then $S = \Bbbk[x_{1},
\dotsc, x_{5}]$ has the $\mathbb{Z}^{2}$-grading defined by 
$\deg(x_{1}) = \deg(x_{2}) = \deg(x_{3}) = \left[ 
\begin{smallmatrix} 1 \\ 0 \end{smallmatrix}
\right]$ and
$\deg(x_{4}) = \deg(x_{5}) = \left[ 
\begin{smallmatrix} 0 \\ 1 \end{smallmatrix}
\right]$.  We consider those multigraded Hilbert polynomials which map
to $3t+1$ under the embedding of $X$ into $\mathbb{P}^{5}$ given by
$\left[
\begin{smallmatrix} 1 \\ 1 \end{smallmatrix}
\right] \in \mathbb{Z}^{2} = \Pic(X)$.
\begin{itemize}
\item $P(t_{1}, t_{2}) = 3t_{1} + 1$: In this case, we need only
consider two decompositions of the multigraded Hilbert polynomial
$P(\bm{t})$, namely $(t_{1} + 1) + (t_{1}) + (t_{1}-1) + 1$ and
$(t_{1} + 1) + (t_{1}) + (t_{1})$.  It follows that the Gotzmann
number is $4$.  Since Proposition~6.10 in \cite{MaclaganSmith1} shows
that $\bm{0} \in \reg(S)$, we deduce that every ideal $I$ with the
given multigraded Hilbert polynomial is $\left[
\begin{smallmatrix} 3 \\ 3 \end{smallmatrix}
\right]$-regular.
\item $P(t_{1}, t_{2}) = 2 t_{1} + t_{2} + 1$: The possible
decompositions are $(t_1+1)+(t_1)+(t_2)$ and
$(t_1+1)+(t_1+1)+(t_2-1)$, so the Gotzmann number is $3$.  
\item $P(t_{1}, t_{2}) = t_{1} + 2 t_{2} + 1$: The only possible
decomposition is $(t_1+1)+(t_2)+(t_2)$, so the Gotzmann number is
again $3$.
\item $P(t_{1}, t_{2}) = 3t_{2} + 1$: There are no $B$-saturated
ideals with this Hilbert polynomial.  Indeed, the first piece of a
decomposition would be $t_2 + 1$ corresponding to a pair $(1, \sigma)$
with $4$, $5 \in \sigma$.  The second pair would have the form
$(x_{i}, \tau)$ for some $i \in \{1,2,3 \}$ which means that the
second piece of the decomposition must again be $t_2 + 1$.  However,
we are left with a polynomial of the form $t_2 - 1$ which is
impossible since we also have $4,5 \in \tau$.
\end{itemize}
\end{example}

%%---------------------------------------------------------------
\section{A New Proof of Gotzmann's Regularity Theorem}
\label{s:standardcase}

By specializing to a standard graded polynomial ring (equivalently to
$\mathbb{P}^{n-1}$), we next show that Theorem~\ref{t:hilbpolybound}
implies Gotzmann's Regularity Theorem.  Throughout this section, $S =
\Bbbk[x_{1}, \dotsc, x_{n}]$ has the $\mathbb{Z}$-grading defined by
$\deg(x_{i}) = 1$ for $1 \leq i \leq n$ and the irrelevant ideal $B =
\langle x_{1}, \dotsc, x_{n} \rangle$.  Gotzmann's Regularity Theorem
gives a bound on the regularity of all $B$-saturated ideals in $S$
with a given Hilbert polynomial $P(t)$.  We first prove that
Gotzmann's bound is the Gotzmann number for $P(t)$ (which justifies
Definition~\ref{d:gotzmannnum}).

\begin{lemma} \label{l:lexsegmentlongest}
If the polynomial $P(t) \in {\mathbb Q}[t]$ can be expressed in the
form
\begin{equation} \label{standarddecomps}
P(t) = \binom{t + q_{1} - u_{1}}{q_{1}} + \binom{t + q_{2} -
u_{2}}{q_{2}} + \dotsb + \binom{t + q_{m} - u_{m}}{q_{m}} \, ,
\end{equation}
where $q_{1} \geq q_{2} \geq \dotsb \geq q_{m} \geq 0$ and $0 \leq
u_{i} \leq i-1$ for $1 \leq i \leq m$ then among all such expressions
the number $m$ is maximized if and only if $u_{i} = i-1$ for all $i$.
\end{lemma}

\begin{proof}
A modification to Algorithm~\ref{a:macaulay} gives a method for
finding all expressions of the form \eqref{standarddecomps}.  Hence,
there is only a finite number of such decompositions, so we may choose
$P(t) = \binom{t + q_{1} - u_{1}}{q_{1}} + \binom{t + q_{2} -
u_{2}}{q_{2}} + \dotsb + \binom{t + q_{m} - u_{m}}{q_{m}}$ to be an
expression of the desired form with a maximal number of summands.
Suppose there is a $i$ such that $u_{i} < i-1$ and let $k$ be the
smallest such $i$.  Using Pascal's identity, we can replace $\binom{t
+ q_{k} - u_{k}}{q_{k}}$ with $\binom{t + q_{k} - (u_{k}+1)}{q_{k}} +
\binom{t + (q_{k}-1) - u_{k}}{q_{k}-1}$.  We claim that by reordering
(if necessary) the binomial coefficients $\binom{t+q_{i} -
u_{i}}{q_{i}}$ with $i > k$ and $\binom{t + (q_{k}-1) -
u_{k}}{q_{k}-1}$, we obtain an expression of the desired form with
$m+1$ summands.  Indeed, the new expression has the desired form
because $u_{k} < k-1$ implies $u_{k+1} \leq k-1$ and the $\binom{t +
(q_{k}-1) - u_{k}}{q_{k}-1}$ term has the same shift with a larger
index.  This longer expression contradicts the maximality of our
choice, however, so we must have $u_{i} = i-1$ for all $i$.
\end{proof}

This lemma allows us to give a new proof of Gotzmann's regularity
theorem.

\begin{theorem}[\cite{gotzmann}] \label{t:gotzmann}
Let $S = \Bbbk[x_{1}, \dotsc, x_{n}]$ be the homogeneous coordinate
ring of $\mathbb{P}^{d}$ and let $B$ be the irrelevant ideal $\langle
x_{1}, \dotsc, x_{n} \rangle$.  If $I$ is an ideal in $S$ and
\begin{equation} \label{gotzdecomp}
P_{S/I}(t) = \binom{t + q_{1}}{q_{1}} + \binom{t + q_{2} - 1} {q_{2}}
+ \dotsb + \binom{t + q_{m} - (m-1)}{q_{m}} \, ,
\end{equation}
where $q_{1} \geq q_{2} \geq \dotsb \geq q_{m} \geq 0$, then $S/(I:
B^{\infty})$ is $(m-1)$-regular.
\end{theorem}

\begin{proof}
By Proposition~\ref{p:usc}, we may assume that $I$ is a $B$-saturated
monomial ideal.  Let $\bigl\{ (\bm{x}^{\bm{u}_{i}}, \sigma_{i}) : 1
\leq i \leq \ell \bigr\}$ be a Stanley filtration for $S/I$ satisfying
the requirements of Proposition~\ref{p:niceSD}.  Since each
$S_{\sigma_{i}}$ is also a standard graded polynomial ring, we know
that each $S_{\sigma_{i}}$ is $0$-regular (see Example~4.2 in
\cite{MaclaganSmith1}).  Remark~\ref{r:totaldegree} implies that $S/I$
is $k$-regular where $k = \max \{ |\bm{u}_{i}| : 1 \leq i \leq m \}$.
We have
\[
P_{S/I}(t) = \binom{t + |\sigma_{1}| - |\bm{u}_{1}|}{|\sigma_{1}|} +
\binom{t + |\sigma_{2}| - |\bm{u}_{2}|}{|\sigma_{2}|} + \dotsb +
\binom{t + |\sigma_{\ell}| - |\bm{u}_{\ell}|}{|\sigma_{\ell}|} \, ,
\]
where $|\sigma_{1}| \geq |\sigma_{2}| \geq \dotsb \geq |\sigma_{\ell}|
\geq 0$ and $0 \leq |\bm{u}_{i}| \leq i-1$ for $1 \leq i \leq \ell$.
Lemma~\ref{l:lexsegmentlongest} shows that $k < \ell \leq m$, which
completes the proof.
\end{proof}

Although not required in our proof of Gotzmann's Regularity Theorem,
the expression \eqref{gotzdecomp} corresponds to a Stanley filtration
of the saturated lexicographic ideal with Hilbert polynomial $P(t)$.
By definition, the $t$th graded component of a lexicographic ideal
$I_{\text{lex}}$ is the $\Bbbk$-vector space spanned by the largest
$H(I_{\text{lex}}, t)$ monomials in lexicographic order.  If we fix an
ordering on the variables $x_{i}$, then Macaulay's description of
Hilbert functions in $S$ (Theorem~4.2.10 in \cite{BH}) shows that
there is a unique $B$-saturated lexicographic ideal associated to
every Hilbert polynomial.

\begin{proposition}
If $P(t) \in {\mathbb Q}[t]$ is a Hilbert polynomial, then the
expression
\begin{equation} \label{gotzmanndecomp}
P(t) = \binom{t + q_{1}}{q_{1}} + \binom{t + q_{2} - 1} {q_{2}} +
\dotsb + \binom{t + q_{m} - (m-1)}{q_{m}} \, ,
\end{equation}
with $q_{1} \geq q_{2} \geq \dotsb \geq q_{m} \geq 0$ comes from a
Stanley filtration for $S/I_{\text{\emph{lex}}}$ where
$I_{\text{\emph{lex}}}$ is the unique $B$-saturated lexicographic ideal
satisfying $P_{S/I_{\text{\emph{lex}}}}(t) = P(t)$.
\end{proposition}

\begin{proof}
From \cite{ReevesStillman}, we know that for every saturated
lexicographic ideal $I_{\text{lex}}$ there is an integer $\ell$
between $0$ and $n-1$ and positive integers $b_{j}$ for $1 \leq j \leq
\ell$ such that
\begin{equation} \label{RSform}
I_{\text{lex}} = \text{{\small $\langle x_{1}, \dotsc, x_{n-\ell-1},
x_{n-\ell}^{b_{1}+1}, x_{n-\ell}^{b_{1}}x_{n-\ell+1}^{b_{2}+1},
\dotsc, x_{n-\ell}^{b_{1}} \dotsb x_{n-3}^{b_{\ell-2}}
x_{n-2}^{b_{\ell-1}+1}, x_{n-\ell}^{b_{1}} \dotsb x_{n-1}^{b_{\ell}}
\rangle$}} \, .
\end{equation}
We use Algorithm~\ref{a:comma-colon} to compute a Stanley filtration
for $S/I$ where $I$ is any ideal of the form given on the right-hand
side of \eqref{RSform}.  In Step~2 of Algorithm~\ref{a:comma-colon},
choose the variable $x_{n - \ell}$; the largest variable dividing the
largest minimal generator of $I$.  It follows that
\begin{align*}
I + \langle x_{n - \ell} \rangle &= \langle x_{1},
\dotsc, x_{n - \ell} \rangle \qquad \text{and} \\ 
\bigl( I : x_{n - \ell} \bigr) &= \langle x_{1}, \dotsc,
x_{n-\ell-1}, x_{n-\ell}^{b_{1}},
x_{n-\ell}^{b_{1}-1}x_{n-\ell+1}^{b_{2}+1}, \dotsc,
x_{n-\ell}^{b_{1}-1} \dotsb x_{n-1}^{b_{\ell}} \rangle \, .
\end{align*}
Hence, the left-hand child of $I$ is prime and corresponds to the pair
$(1, \{ n - \ell + 1, \dotsc, n \})$.  On the other hand, the
right-hand child is another ideal of the form given on the right-hand
side of \eqref{RSform}. Iterating this process, we obtain a Stanley
filtration of $S/I$:
\begin{equation} \label{sdforlex}
\bigcup_{j = 1}^{\ell} \bigcup_{i = 0}^{b_{j}} \bigl\{
(x_{n-\ell}^{b_{1}} \dotsb x_{n-\ell+j-1}^{i}, \{ n-\ell+j, \dotsc, n
\}) \bigr\} \, ,
\end{equation}
and one easily verifies that \eqref{sdforlex} yields an expression of
the form \eqref{gotzmanndecomp}.  
\end{proof}

Since the number of pairs in the Stanley filtration \eqref{sdforlex}
equals the maximum total degree of a minimal generator of the
saturated lexicographic ideal $I_{\text{lex}}$, it follows that
Gotzmann's regularity theorem is sharp. This establishes the
well-known result that the lexicographic ideal has the worst
regularity among all $B$-saturated ideals with the same Hilbert
polynomial.

%%---------------------------------------------------------------
\section{Multigraded Hilbert schemes} \label{s:hilbschemes}

The aim of this section is to construct a space $\Hilb_{X}^{P}$ that
parameterizes all subschemes of $X$ with a given multigraded Hilbert
polynomial $P \in \mathbb{Q}[t_{1}, \dotsc, t_{r}]$.  This generalizes
the original Hilbert scheme, introduced in \cite{G}, which
parameterizes subschemes of projective space.  Like all parameter
spaces, $\Hilb_{X}^{P}$ allows one to study the natural adjacency
relationships between subschemes.  This larger class of Hilbert
schemes also includes many more manageably sized examples.  By
analyzing these small spaces, especially those which are accessible to
computational experimentation, we expect to gain new insights into
Hilbert schemes.

Before discussing our construction, we provide a simple example.

\begin{example}
It is well-known that the lines on the nonsingular quadratic surface
$X \cong \mathbb{P}^{1} \times \mathbb{P}^{1}$ contained in
$\mathbb{P}^{3}$ belong to two families.  In fact, each of these
families is precisely a multigraded Hilbert scheme.  Specifically, the
closed subscheme of $\Hilb_{\mathbb{P}^{3}}^{t+1}$ parameterizing
subschemes of $\mathbb{P}^{3}$ with Hilbert polynomial $t+1$ lying on
$X$ is the disjoint union $\Hilb_{X}^{t_{1}+1} \amalg
\Hilb_{X}^{t_{2}+1} \cong \mathbb{P}^{1} \amalg \mathbb{P}^{1}$.
\end{example}

We construct the space $\Hilb_{X}^{P}$ by proving that the appropriate
functor is represented by a projective scheme.  Define the functor
$\Hilb_{X}^{P}$ that sends the category of commutative rings over
$\Bbbk$ to the category sets as follows: given a commutative ring $R$
over $\Bbbk$, $\Hilb_{X}^{P}(R)$ is the set of families of subschemes
$Y \subseteq X \times_{\Bbbk} \Spec(R)$ over $\Spec(R)$ whose sheaf of
ideals has the specified Hilbert polynomial $P$.  To prove that
$\Hilb_{X}^{P}$ is representable, we build on the methods used in
\cite{HS}; see \S6.1 for the explicit reference to our setting.

To begin, we recall the Hilbert functor $H_{S_{\scr{D}}}^{h}$ from
\cite{HS}.  For a subset $\scr{D} \subset \mathbb{Z}^{r}$, we write
$S_{\scr{D}}$ for the graded $\Bbbk$-vector space
$\bigoplus\nolimits_{\bm{p} \in \scr{D}} S_{\bm{p}}$ and $F_{\scr{D}}
= \bigcup\nolimits_{\bm{p}, \bm{k} \in \scr{D}} F_{\bm{p}, \bm{k}}$
denotes a collection of maps from $S_{\bm{p}}$ to $S_{\bm{k}}$.  More
precisely, $F_{\bm{p}, \bm{k}}$ consists of the multiplication maps
arising from the monomials in $S_{\bm{k} - \bm{p}}$.  For a
commutative ring $R$ over $\Bbbk$, let $R \otimes S_{\scr{D}}$ be the
graded $R$-module $\bigoplus_{\bm{p} \in \scr{D}} R \otimes_{\Bbbk}
S_{\bm{p}}$ with operators $F_{\bm{p}, \bm{k}}^{R} = (1_{R}
\otimes_{\Bbbk} -)(F_{\bm{p}, \bm{k}})$.  A homogeneous submodule $L =
\bigoplus_{\bm{p} \in \scr{D}} L_{\bm{p}} \subseteq R \otimes
S_{\scr{D}}$ is an $F$-submodule if it satisfies $F_{\bm{p},
\bm{k}}^{R} (L_{\bm{p}}) \subseteq L_{\bm{k}}$ for all $\bm{p}, \bm{k}
\in \scr{D}$.  Given a function $h \colon \scr{D} \longrightarrow
\mathbb{N}$, let $H_{S_{\scr{D}}}^{h}(R)$ be the set of $F$-submodules
$L \subseteq R \otimes S_{\scr{D}}$ such that $(R \otimes_{\Bbbk}
S_{\bm{p}})/ L_{\bm{p}}$ is a locally free $R$-module of rank
$h(\bm{p})$ for each $\bm{p} \in \scr{D}$.  If $\psi \colon R
\longrightarrow R'$ is a homomorphism, then local freeness implies
that $L' = R' \otimes_{R} L$ is an $F$-submodule of $R' \otimes
S_{\scr{D}}$ and $(R' \otimes_{\Bbbk} S_{\bm{p}})/ L_{\bm{p}}'$ is a
locally free $R'$-module of rank $h(\bm{p})$ for each $\bm{p} \in
\scr{D}$.  Defining $H_{S_{\scr{D}}}^{h}(\psi) \colon
H_{S_{\scr{D}}}^{h}(R) \longrightarrow H_{S_{\scr{D}}}^{h}(R')$ to be
the map sending $L$ to $L'$ makes $H_{S_{\scr{D}}}^{h}$ into a functor
from the category of commutative rings over $\Bbbk$ to the category of
sets.

When the function $h \colon \scr{D} \longrightarrow \mathbb{N}$ is
defined by evaluating a polynomial $P$ at points in $\scr{D}$, we
simply write $H_{S_{\scr{D}}}^{P}$.  By relating the functors
$\Hilb_{X}^{P}$ and $H_{S_{\scr{D}}}^{P}$, we show that
$\Hilb_{X}^{P}$ is representable.

\begin{theorem}
If $P \in \mathbb{Q}[t_{1}, \dotsc, t_{r}]$ is a Hilbert polynomial,
then the functor $\Hilb_{X}^{P}$ is represented by a projective scheme
over $\Bbbk$.  In fact, there is a finite subset $\scr{D} \subset
\mathbb{Z}^{r}$ which produces a canonical closed embedding from
$\Hilb_{X}^{P}$ into $H_{S_{\scr{D}}}^{P}$.
\end{theorem}

\begin{proof}
If $R$ is a commutative ring over $\Bbbk$, then \cite{Cox} shows that
each ideal sheaf in $\Hilb_{X}^{P}(R)$ corresponds to unique
$B$-saturated ideal $I$ in the ring $S \otimes_{\Bbbk} R = R[x_{1},
\dotsc, x_{n}]$.  Using Theorem~\ref{t:hilbpolybound}, we can choose a
$\bm{k} \in \scr{K}$ for which every such $I$ is $\bm{k}$-regular.
Lemma~6.8 in \cite{MaclaganSmith1} states that the truncation
$I|_{\bm{k} + \scr{K}} := S \cdot \bigl( \bigoplus_{\bm{p} \in \bm{k}
+ \scr{K}} I_{\bm{p}} \bigr)$ corresponds to the same ideal sheaf on
$X$ as $I$ does.  This bijection between sheaves of ideals on $X$ and
truncations of ideals in $S$ gives a natural transformation between
$\Hilb_{X}^{P}$ and $H_{S_{\bm{k} + \scr{K}}}^{P}$.

In \S6.1 of \cite{HS} Haiman and Sturmfels claim that there exists a
finite set $\scr{D} \subset \bm{k} + \scr{K}$ satisfying
\begin{equation} \label{supportive}
\text{
\begin{minipage}{13cm}
for every extension field $K$ of $\Bbbk$ and every $L_{\scr{D}} \in
H_{S_{\scr{D}}}^{P}(K)$, if $L'$ denotes the $F$-submodule of $K
\otimes S_{\scr{D}}$ generated by $L_{\scr{D}}$ then $\dim (K
\otimes_{\Bbbk} S_{\bm{t}}) / L_{\bm{t}}' \leq P(\bm{t})$ for all
$\bm{t} \in \bm{k} + \scr{K}$. 
\end{minipage}}
\end{equation}
For such a finite set $\scr{D} \subset \bm{k} + \scr{K}$, Theorem~2.3
in \cite{HS} produces a closed embedding $\Hilb_{X}^{P} = H_{S_{\bm{k}
+ \scr{K}}}^{P} \longrightarrow H_{S_{\scr{D}}}^{P}$.  Since
Theorem~2.2 and Remark~2.5 in \cite{HS} prove that
$H_{S_{\scr{D}}}^{P}$ is represented by a closed subscheme of a
Grassmann scheme, this completes the proof.
\end{proof}

To give explicit equations for $\Hilb_{X}^{P}$, we need an effective
description of both the set $\scr{D}$ and the equations defining the
closed subscheme of $H_{S_{\scr{D}}}^{P}$.  The following algorithm,
essentially a constructive version of Proposition~3.2 in \cite{HS},
produces the set $\scr{D}$.

\begin{algorithm}
Given a Hilbert polynomial $P \in \mathbb{Q}[t_{1}, \dotsc, t_{r}]$,
this algorithm returns a finite subset $\scr{D}$ satisfying
\eqref{supportive}.
\begin{enumerate}
\item (Initialize) Set $\scr{D}$ equal to $\{\bm{k}\}$, where $\bm{k}$
is a bound on the regularity of all ideals with Hilbert polynomial $P$
obtained from Theorem~\ref{t:hilbpolybound}.
\item (Create ideals) Construct the set $\textsf{Ideals}$ of all
monomial ideals $I$ generated in degree $\scr{D}$ such that $H(S/I,
\bm{t}) = P(\bm{t})$ for all $\bm{t} \in \scr{D}$.  Since there are
only a finite number of monomials with degrees in $\scr{D}$, this is a
finite set.
\item (Finished?) If every ideal $I$ in $\textsf{Ideals}$ satisfies
$P_{S/I}(\bm{t}) = P(\bm{t})$ then return $\scr{D}$.  Otherwise, for
every ideal in $\textsf{Ideals}$ find a $\bm{t} \in \bm{k} + \scr{K}$
such that $H(S/I, \bm{t}) \neq P(\bm{t})$.  Add each of these points
to $\scr{D}$ and return to Step~2.  One choice of such points is to
use the maximum degree of a monomial with degree in $\scr{D}$ to bound
the maximum size of any $|\bm{u}_i|$, and thus of any $A\bm{u}_i$,
occurring in a Stanley filtration of the appropriate form.  This gives
a bound $\bm{c}$ on the regularity of all ideals generated in
$\scr{D}$, and so we can add the point $\bm{c}$, together with
$\binom{n}{d}$ sufficiently general points in $\bm{c} + \scr{K}$ to
$\scr{D}$.  Evaluating $H(S/I, \bm{t})$ at these points also lets us
check whether $P_{S/I}(\bm{t}) = P(\bm{t})$.
\end{enumerate}
\end{algorithm}

\begin{proof}[Proof of Correctness]
The proof of Proposition 3.2 of \cite{HS} establishes that this
algorithm terminates.  It remains to show that the output satisfies
\eqref{supportive}.  By construction, every ideal $I$ in
$\textsf{Ideals}$ has Hilbert polynomial $P$.  Step~1 guarantees that
the saturation $\overline{I} = (I : B^{\infty})$ has Hilbert
polynomial $P$ and is $\bm{k}$-regular.  Theorem~5.4 in
\cite{MaclaganSmith1} implies that $\overline{I}|_{\bm{k}+\scr{K}}$ is
generated in degree $\bm{k}$.  Since $H(S/I,\bm{k}) =
H(S/\overline{I}, \bm{k}) = P(\bm{k})$, we have $I_{\bm{k}} =
\overline{I}_{\bm{k}}$.  Because $I \subseteq \overline{I}$, it
follows that $I|_{\bm{k} + \scr{K}} = \overline{I}|_{\bm{k} +
\scr{K}}$.  Applying Corollary~\ref{c:regmeanspoly}, we see that
$H(S/\overline{I}, \bm{t}) = P(\bm{t})$ for all $\bm{t} \in
\bm{k}+\scr{K}$.  We conclude that \eqref{supportive} holds.

We finish by explaining why in the Step~3 it suffices to choose
$\binom{n}{d}$ sufficiently general points in $\bm{c}+\scr{K}$ to add
to $\scr{D}$.  By construction all ideals generated in $\scr{D}$ agree
with their Hilbert polynomial on $\bm{c} + \scr{K}$.  Since
$P(\bm{t})$ is a polynomial of degree at most $d$ in $r$ variables, it
has at most $\binom{d+r}{d}$ terms.  If the Hilbert function of an
ideal $I$ generated in the degrees in $\scr{D}$ agrees with
$P(\bm{t})$ for $\binom{n}{d}$ sufficiently general points in
$\bm{c}+\scr{K}$, then it must have Hilbert polynomial $P$.
\end{proof}

A multigraded version of Gotzmann's Persistence Theorem would lead to
an effective description of the equations defining the relevant closed
subscheme of $H_{S_{\scr{D}}}^{P}$.  This is the central open problem
in this area.

%%--REFERENCES----------------------------------------------
%\bibliographystyle{amsalpha}
%\bibliography{toricRegularity}

\begin{thebibliography}{LHTY}

\bibitem[ACD]{ACdN}
Annetta Aramova, Kristina Crona, and Emanuela De~Negri, 
  \emph{Bigeneric initial ideals, diagonal subalgebras and bigraded 
  {H}ilbert functions}, Journal of Pure and Applied Algebra 
  \textbf{150} (2000), no.~3, 215--235.

\bibitem[AK]{altmanKleiman}
Allen~B. Altman and Steven~L. Kleiman, \emph{Compactifying the {P}icard
  scheme}, Adv. in Math. \textbf{35} (1980), no.~1, 50--112.

\bibitem[Ape]{apel}
Joachim Apel, \emph{On a conjecture of {R}. {P}. {S}tanley; part {II} -
  {Q}uotients modulo monomial ideals}, MSRI Preprint \#2001-009, 2001.

\bibitem[Bat]{Batyrev}
Victor~V. Batyrev, \emph{On the classification of toric {F}ano 
  {$4$}-folds}, Algebraic geometry~9, Journal of Mathematical Sciences
  (New York) \textbf{94} (1999), no.~1, 1021--1050.

\bibitem[BH]{BH}
Winfried Bruns and J{\"u}rgen Herzog, \emph{Cohen-{M}acaulay rings}, 
  Cambridge Studies in Advanced Mathematics, vol.~39, Cambridge 
  University Press, Cambridge, 1993.

\bibitem[BS]{BrodmannSharp}
Markus~P. Brodmann and Rodney~Y. Sharp, \emph{Local cohomology: an 
  algebraic introduction with geometric applications}, Cambridge 
  Studies in Advanced Mathematics, vol.~60, Cambridge University 
  Press, Cambridge, 1998.

\bibitem[Cox1]{Cox}
David~A. Cox, \emph{The homogeneous coordinate ring of a toric 
  variety}, Journal of Algebraic Geometry \textbf{4} (1995), no.~1, 
  17--50.

\bibitem[Cox2]{cox2}
\bysame, \emph{Recent developments in toric geometry}, Algebraic
  geometry---Santa Cruz 1995, Proc. Sympos. Pure Math., vol.~62, Amer.
  Math. Soc., Providence, RI, 1997, pp.~389--436.

\bibitem[EL]{EL}
Lawrence Ein and Robert Lazarsfeld, \emph{Syzygies and {K}oszul 
  cohomology of smooth projective varieties of arbitrary dimension}, 
  Inventiones Mathematicae \textbf{111} (1993), no.~1, 51--67.

\bibitem[Ful]{fulton}
William Fulton, \emph{Introduction to toric varieties}, Annals of 
  Mathematics Studies~131, Princeton University Press, Princeton, NJ, 
  1993, The William H. Roever Lectures in Geometry.

\bibitem[GLP]{GLP}
Laurent Gruson, Robert Lazarsfeld, and Christian Peskine, \emph{On a 
  theorem of {C}astelnuovo, and the equations defining space curves}, 
  Inventiones Mathematicae \textbf{72} (1983), no.~3, 491--506.

\bibitem[Got]{gotzmann}
Gerd Gotzmann, \emph{Eine {B}edingung f{\"{u}}r die {F}lachheit und das
  {H}ilbertpolynom eines graduierten {R}inges}, Mathematische 
  Zeitschrift \textbf{158} (1978), no.~1, 61--70.

\bibitem[Gre1]{GreenGot}
Mark Green, \emph{Restrictions of linear series to hyperplanes, and 
  some results of {M}acaulay and {G}otzmann}, Algebraic curves and 
  projective geometry (Trento, 1988), Lecture Notes in Math., vol. 
  1389, Springer, Berlin, 1989, pp.~76--86.

\bibitem[Gre2]{Green}
Mark~L. Green, \emph{Generic initial ideals}, Six lectures on 
  commutative algebra (Bellaterra, 1996), Birkh\"auser, Basel, 1998, 
  pp.~119--186.

\bibitem[Gro]{G}
Alexander Grothendieck, \emph{Techniques de construction et 
  th\'eor\`emes d'existence en g\'eom\'etrie alg\'ebrique. {I}{V}. 
  {L}es sch\'emas de {H}ilbert}, S\'eminaire Bourbaki, Vol.\ 6, Soc. 
  Math. France, Paris, 1995, pp.~Exp.\ No.\ 221, 249--276.

\bibitem[M2]{M2}
Daniel~R. Grayson and Michael~E. Stillman, \emph{Macaulay 2, a 
  software system for research in algebraic geometry}, Available at
  \texttt{http://www.math.uiuc.edu/Macaulay2}.

\bibitem[HS]{HS}
Mark Haiman and Bernd Sturmfels, \emph{Multigraded {H}ilbert schemes},
  \texttt{arXiv:math.AG/0201271}.

\bibitem[HTr]{HoangTrung}
Nguyen~Duc Hoang and Ngo~Viet Trung, \emph{{Hilbert polynomials of 
  non-standard bigraded algebras}}, \texttt{arXiv:math.AC/0211181}.

\bibitem[HTh]{HoTh1}
Serkan Ho{\c{s}}ten and Rekha~R. Thomas, \emph{Standard pairs and 
  group relaxations in integer programming}, Journal of Pure and 
  Applied Algebra \textbf{139} (1999), no.~1-3, 133--157, Effective 
  methods in algebraic geometry (Saint-Malo, 1998).

\bibitem[Kle]{Kleiman}
Steven~L. Kleiman, \emph{Toward a numerical theory of ampleness}, 
  Annals of Mathematics. Second Series \textbf{84} (1966), 293--344.

\bibitem[LHTY]{LattE}
Jes\'{u}s A.~De Loera, Raymond Hemmecke, Jeremiah Tauzer, and Ruriko
  Yoshida, \newline
  \emph{Effective lattice point counting in rational convex
  polytopes}, \texttt{LattE} code available at \newline
  \texttt{http://www.math.ucdavis.edu/\textasciitilde latte}.

\bibitem[Mac]{antichains}
Diane Maclagan, \emph{Antichains of monomial ideals are finite}, 
  Proceedings of the American Mathematical Society \textbf{129} 
  (2001), no.~6, 1609--1615.

\bibitem[MS]{MaclaganSmith1}
Diane Maclagan and Gregory~G. Smith, \emph{Multigraded 
  {C}astelnuovo-{M}umford regularity}, 2003.

\bibitem[Mum]{Mumford}
David Mumford, \emph{Lectures on curves on an algebraic surface}, 
  Princeton University Press, Princeton, N.J., 1966.

\bibitem[RS]{ReevesStillman}
Alyson Reeves and Mike Stillman, \emph{Smoothness of the 
  lexicographic point}, Journal of Algebraic Geometry \textbf{6} 
  (1997), no.~2, 235--246.

\bibitem[Sch]{schrijver}
Alexander Schrijver, \emph{Theory of linear and integer programming}, 
  John Wiley \& Sons Ltd., Chichester, 1986, A Wiley-Interscience 
  Publication.

\bibitem[Sim]{clean}
Robert~Samuel Simon, \emph{Combinatorial properties of ``cleanness''},
  Journal of Algebra \textbf{167} (1994), no.~2, 361--388.

\bibitem[Sta]{StanleySD}
Richard~P. Stanley, \emph{Linear {D}iophantine equations and local 
  cohomology}, Inventiones Mathematicae \textbf{68} (1982), no.~2, 
  175--193.

\bibitem[St1]{SturmfelsVPF}
Bernd Sturmfels, \emph{On vector partition functions}, Journal of 
  Combinatorial Theory. Series A \textbf{72} (1995), no.~2, 302--309.

\bibitem[St2]{sturmfels}
\bysame, \emph{Gr\"obner bases and convex polytopes}, American 
  Mathematical Society, Providence, RI, 1996.

\bibitem[STV]{STV}
Bernd Sturmfels, Ng{\^o}~Vi{\^e}t Trung, and Wolfgang Vogel, 
  \emph{Bounds on degrees of projective schemes}, Mathematische 
  Annalen \textbf{302} (1995),
  no.~3, 417--432.

\bibitem[SW]{SturmfelsWhite}
Bernd Sturmfels and Neil White, \emph{Computing combinatorial 
  decompositions of rings}, Combinatorica \textbf{11} (1991), no.~3, 
  275--293.

\bibitem[Vie]{Viehweg}
Eckart Viehweg, \emph{Quasi-projective moduli for polarized manifolds},
  Springer-Verlag, Berlin, 1995.

\bibitem[Zie]{Ziegler}
G{\"u}nter~M. Ziegler, \emph{Lectures on polytopes}, Graduate Texts in
  Mathematics, vol. 152, Springer-Verlag, New York, 1995.

\end{thebibliography}

\def\cprime{$'$}
\providecommand{\bysame}{\leavevmode\hbox to3em{\hrulefill}\thinspace}

\end{document}